\documentclass[leqno]{siamart171218}
\usepackage{amssymb}
\usepackage{graphicx,xcolor} 
\usepackage{xspace}
\usepackage{multirow}
\usepackage{bbm}
\usepackage[numbers]{natbib}
\let\cite=\citet
\usepackage{enumerate}
\usepackage{bm} 
\newif\ifOUT\OUTfalse
\ifOUT
\newcommand{\Question}[1]{}
\newcommand{\done}[1]{{}}
\newcommand{\afaire}[1]{{}}

\else
\newcommand{\Question}[1]{{\marginpar{\color{blue}#1}}}
\newcommand{\done}[1]{{\marginpar{\color{red}$\heartsuit\heartsuit$#1}}}
\newcommand{\afaire}[1]{{\marginpar{\color{blue}$\spadesuit\spadesuit$#1}}}

\definecolor{verde}{rgb}{0.2, 0.5, 0.3}

\fi

\numberwithin{theorem}{section}
\numberwithin{equation}{section}

\theoremstyle{plain} 
\theoremheaderfont{\normalfont\scshape}
\theorembodyfont{\normalfont\itshape} 
\theoremseparator{.}
\theoremsymbol{\hfill\ensuremath{\square}}
\newtheorem{Th}[theorem]{Theorem}

\newtheorem{Lem}[theorem]{Lemma}

{\theoremheaderfont{\normalfont\itshape}
\theorembodyfont{\normalfont}
\newtheorem{Rem}[theorem]{Remark}
\newtheorem{Exp}[theorem]{Example}}

\def\bi{\begin{itemize}}
\def\ei{\end{itemize}}
\def\berom{\begin{enumerate}[{\rm(i)}]}
\def\eerom{\end{enumerate}}

\newcounter{itemrem}

\def\eRem{\qed\end{Rem}}
\makeatletter
\def\bRem{\@ifnextchar[{\@remwithtitle}{\@remwithouttitle}}
\def\@remwithtitle[#1]{\begin{Rem}[#1]\setcounter{itemrem}{0}}%
\def\@remwithouttitle{\begin{Rem}\setcounter{itemrem}{0}}%
\makeatother

\newcounter{itemexp}

\def\eExp{\qed\end{Exp}}
\makeatletter
\def\bExp{\@ifnextchar[{\@expwithtitle}{\@expwithouttitle}}
\def\@expwithtitle[#1]{\begin{Exp}[#1]\setcounter{itemexp}{0}}%
\def\@expwithouttitle{\begin{Exp}\setcounter{itemexp}{0}}%
\makeatother

\def\bproof{\begin{proof}}
\def\eproof{\qed\end{proof}}

\ifOUT

\else

\fi

\def\bExo{\begin{Exo}}
\def\eExo{\end{Exo}}

\def\wc{{\widehat c}}

\def\wF{{\widehat F}}

\def\wK{{\widehat K}}

\def\wP{{\widehat P}}

\def\wbx{{\widehat \bx}}

\def\wbP{{\widehat \bP}}

\def\wSigma{{\widehat {\mit\Sigma}}}

\def\tp{{\widetilde p}}

\newcommand{\bfjlg}{\bm}         
\newcommand{\boldsymboljlg}{\boldsymbol} 

\newcommand{\bb}{{\bfjlg b}}

\newcommand{\bn}{{\bfjlg n}}

\newcommand{\bv}{{\bfjlg v}}
\newcommand{\bw}{{\bfjlg w}}
\newcommand{\bx}{{\bfjlg x}}

\newcommand{\bC}{{\bfjlg C}}

\newcommand{\bH}{{\bfjlg H}}

\newcommand{\bL}{{\bfjlg L}}

\newcommand{\bP}{{\bfjlg P}}

\newcommand{\bS}{{\bfjlg S}}
\newcommand{\bT}{{\bfjlg T}}

\newcommand{\bV}{{\bfjlg V}}
\newcommand{\bW}{{\bfjlg W}}

\newcommand{\bY}{{\bfjlg Y}}

\newcommand{\bbeta}{{\boldsymboljlg \beta}}

\newcommand{\bsigma}{{\boldsymboljlg \sigma}}

\newcommand{\btau}{{\boldsymboljlg \tau}}
\newcommand{\bphi}{{\boldsymboljlg \phi}}
\newcommand{\bpsi}{{\boldsymboljlg \psi}}

\newcommand{\bzero}{\bfjlg{0}}

\newcommand{\calE}{{\mathcal E}}
\newcommand{\calF}{{\mathcal F}}

\newcommand{\calH}{{\mathcal H}}
\newcommand{\calI}{{\mathcal I}}

\newcommand{\calK}{{\mathcal K}}
\newcommand{\calL}{{\mathcal L}}

\newcommand{\calT}{{\mathcal T}}

\newcommand{\polN}{{\mathbb N}}

\newcommand{\polP}{{\mathbb P}}

\newcommand{\upb}{^{\mathrm{b}}}
\newcommand{\upc}{^{\mathrm{c}}}
\newcommand{\upd}{^{\mathrm{d}}}

\newcommand{\upg}{^{\mathrm{g}}}

\newcommand{\upgav}{^{\mathrm{g,av}}}

\newcommand\opRec{{\mathsf R}}

\newcommand{\Vshs}{V_{\sharp}}

\def\Real{{\mathbb R}}

\def\form{\ell}

\def\inter{{\cal I}}

\renewcommand{\det}{\operatorname{det}}
\renewcommand{\dim}{\operatorname{dim}}

\def\div{\nabla{\cdot}}

\def\lb{[\![}
\def\rb{]\!]}

\newcommand{\upint}{^\circ}
\newcommand{\upbnd}{^\partial}

\def\hinH{{h\in\calH}}
\def\subhinH{_\hinH}
\def\famTh{(\calT_h)\subhinH}

\def\calFh{\calF_h}
\def\calFK{\calF_K}
\def\calFhi{\calFh\upint}
\def\calFhb{\calFh\upbnd}

\def\wKPS{(\wK,\wP,\wSigma)}
\def\wbKPS{(\wK,\wbP,\wSigma)}

\def\uCR{^{\textsc{cr}}}

\newenvironment{Ventry}[1]%
{\begin{list}{}{%
\settowidth{\labelwidth}{{\rm (#1)}}%
\setlength{\leftmargin}{\labelwidth+\labelsep}}}
{\end{list}}

\newcommand{\avg}[1]{\{#1\}}
\newcommand{\jumpsn}[1]{|#1|_{\mathrm{J}}}
\newcommand{\bndsn}[1]{|#1|_{\partial}}

\newcommand{\Dom}{D}
\newcommand{\Jac}{\mathbb{J}}

\newcommand{\difT}{\mathbbm{d}}

\newcommand{\calThbb}{\overline{\calT}_h^{\partial}}

\newcommand{\intset}[2]{\{#1\hskip.05em\relax{:}\hskip.1em\relax#2\}}

\newcommand{\mes}[1]{|#1|}

\newcommand{\term}{\mathfrak{T}}

\newcommand{\opSt}{{\mathsf S}}

\newcommand{\loS}{_{\mathrm{\scriptscriptstyle S}}}

\def\DIV{\nabla{\cdot}}      
\def\ROT{\nabla{\times}} 
\def\GRAD{\nabla}

\def\SCAL{{\cdot}}

\def\CROSS{{\times}}

\def\dif{\,\mathrm{d}}

\def\tr{^{\sf T}}

\def\front{{\partial\Dom}}

\def\Ldeuxd{{\bL^2(\Dom)}}

\def\Ldeux{{L^2(\Dom)}}

\def\Hunz{{H^1_0(\Dom)}}
\def\Hmun{{H^{-1}(\Dom)}}

\def\Hun{{H^{1}(\Dom)}}

\def\Hrot{{{\bH(\text{\rm curl};\Dom)}}}

\def\jump#1{\lb{#1}\rb}

\newcommand{\mapK}{{\psi_K}}

\newcommand{\mapKg}{\psi_K\upg}
\newcommand{\mapFg}{\psi_F\upg}
\newcommand{\mapKc}{\bpsi_K\upc}

\newcommand{\ie}{i.e.,\@\xspace}
\newcommand{\eg}{e.g.,\@\xspace}

\newcommand{\sth}{s.t.\@\xspace}
\renewcommand{\ae}{a.e.\@\xspace}
\newcommand{\wrt}{w.r.t.\@\xspace}

\newcommand{\bVentry}[1]{\begin{Ventry}{#1}}
\newcommand{\eVentry}{\end{Ventry}}

\catcode`\@=11
\newdimen\linespacing
\linespacing=\baselineskip
\newcommand\subsubsectionmodif{\@startsection{subsubsection}{3}%
  {\z@}{.5\linespacing\@plus.7\linespacing}{-.5em}{\bfseries}*}

\newlength\pageone
\newlength\pagetwo
\newlength\interspace

\newlength\retraitspace
\newlength\calculspace

\newcommand{\tq}{{\;|\;}} 
\newcommand{\st}{\tq}
\newcommand{\bset}{\{} \newcommand{\eset}{\}}



\makeatletter
\newcommand{\manuallabel}[2]{{\textup{(#2)}\def\@currentlabel{#2}\label{#1}}}
\makeatother

\newcounter{subeq}

\begin{document}
\newcommand\footnotemarkfromtitle[1]{%
\renewcommand{\thefootnote}{\fnsymbol{footnote}}%
\footnotemark[#1]%
\renewcommand{\thefootnote}{\arabic{footnote}}}

\title{Quasi-optimal nonconforming approximation of elliptic PDEs with contrasted coefficients and minimal regularity\footnotemark[1]}

\author{Alexandre Ern\footnotemark[2] \and Jean-Luc Guermond\footnotemark[3],
\date{Draft version \today}}

\maketitle

\renewcommand{\thefootnote}{\fnsymbol{footnote}} \footnotetext[1]{
  This material is based upon work supported in part by the National
  Science Foundation grants DMS-1619892, DMS-1620058, by the Air Force
  Office of Scientific Research, USAF, under grant/contract number
  FA9550-18-1-0397, and by the Army Research Office under
  grant/contract number W911NF-15-1-0517.}
\footnotetext[2]{Department of Mathematics, Texas A\&M University 3368
  TAMU, College Station, TX 77843, USA.}
\footnotetext[2]{Universit\'e Paris-Est, CERMICS (ENPC),
  77455 Marne-la-Vall\'ee cedex 2, France and INRIA Paris, 2 rue Simone Iff, 75589 Paris, France.}
\renewcommand{\thefootnote}{\arabic{footnote}}

\begin{abstract} 
  In this paper we investigate the approximation of a diffusion model
  problem with contrasted diffusivity and the error analysis
  of various nonconforming approximation methods. The essential
  difficulty is that the Sobolev smoothness index of the exact
  solution may be just barely larger than one.  The lack of smoothness
  is handled by giving a weak meaning to the normal derivative of the
  exact solution at the mesh faces. The error estimates are
  robust with respect to the diffusivity contrast. We briefly show how
  the analysis can be extended to the Maxwell's equations.
\end{abstract}



\pagestyle{myheadings} \thispagestyle{plain} \markboth{A. ERN, J.-L. GUERMOND}{Contrasted diffusion}

\begin{center}
\textit{This article is dedicated to the memory of Christine Bernardi.}
\end{center}

\section{Introduction} \label{Sec:introduction} 

The objective of the present paper is to revisit and unify the error
analysis of various nonconforming approximation techniques applied to
a diffusion model problem with contrasted diffusivity. We also briefly
show how to extend the analysis to Maxwell's equations.

\subsection{Content of the paper}

The nonconforming techniques we have in mind are Crouzeix--Raviart
finite elements \citep{CR73}, Nitsche's boundary penalty method
\citep{Nitsche_1971}, the interior penalty discontinuous Galerkin
(IPDG) method \citep{Arnol:82}, and the hybrid high-order (HHO)
methods \citep{DiPEr:15,DiPEL:14} which are closely
related to hybridizable discontinuous Galerkin methods
\citep{CoDPE:16}.  The main difficulty in the error analysis is that
owing to the contrast in the diffusivity, the Sobolev smoothness index
of the exact solution is barely larger than one.  This makes the
estimation of the consistency error incurred by nonconforming
approximation techniques particularly challenging since the normal
derivative of the solution at the mesh faces is not integrable and it
is thus not straightforward to give a reasonable meaning to this quantity on each
mesh face independently.

The main goal of the present paper is to establish quasi-optimal
error estimates by using a mesh-dependent norm that remains bounded as
long as the exact solution has a Sobolev smoothness index strictly larger 
than one.
By quasi-optimality, we mean that the approximation error measured in
the augmented norm is bounded, up to a generic constant, by the best
approximation error of the exact solution measured in the same
augmented norm by members of the discrete trial space.  A key point in
the analysis is that the above generic constant is independent of the
diffusivity contrast. We emphasize that quasi-optimal error estimates
are more informative than the more traditional asymptotic error
estimates, which bound the approximation error by terms that optimally
decay with the mesh size. Indeed, the former estimates cover
the whole computational range whereas the latter estimates only
  cover the asymptotic range.  One key novelty herein is the
introduction of a weighted bilinear form that accounts for the default
of consistency in all the cases (see \eqref{eq:def_n_sharp}).

The paper is organized as follows. The model problem under
consideration and the discrete setting are introduced in
\S\ref{Sec:Preliminaries}.  The weighted bilinear form mentioned above
which accounts for the consistency default at the mesh interfaces and
boundary faces is defined in \S\ref{Sec:nsharp}. The key results in
this section are Lemma~\ref{lem:identity_n_sharp_dif} and
Lemma~\ref{lem:bnd_n_sharp}.  We collect in \S\ref{Sec:applications}
the error analyses of the approximation of the model problem with the
Crouzeix--Raviart approximation, Nitsche's boundary penalty method,
the IPDG approximation, and the HHO approximation.
To avoid invoking Strang's second Lemma, we introduce in
\S\ref{Sec:abstract_error} a linear form $\delta_h$ that measures
consistency but does not need the exact solution to be inserted into
the arguments of the discrete bilinear form at hand.  The weighted
bilinear form \eqref{eq:def_n_sharp} turns out to an essential tool to
deduce robust estimates of the norm of the consistency form $\delta_h$
for all the nonconforming methods considered.  One originality of this
paper is that all the error estimates provided in
\S\ref{Sec:applications} involve constants that are uniform with
respect to the diffusivity contrast. Another salient feature is
that the source term is assumed to
be only in $L^q(\Dom)$, where $q$ is such that $L^{q}(\Dom)$ is continuously
embedded in $\Hmun:=(\Hunz)'$; specifically, this means that 
$q>2_*:=\frac{2d}{2+d}\ge 1$ 
(here, $d\ge2$ is the space dimension).

\subsection{Literature overview}
Let us put our work in perspective with the literature.  
Perhaps a bit surprisingly, error estimates for nonconforming
approximation methods are rarely presented in a quasi-optimal form
in the literature. A key step
toward achieving quasi-optimal error estimates has been achieved  in
\cite{VeeZa:18,VeeZa:18b}. Therein, the approximation error and the
best-approximation error are both measured using the energy norm and
the source term is assumed to be just in the dual space
$H^{-1}(\Dom)$. However, at the time of this
writing, 
this setting does not yet cover robust estimates \wrt the diffusivity
contrast. In the present work, we proceed somewhat differently to
obtain robust quasi-optimal error estimates.  This is done at the
following price: (i) We invoke augmented norms, which are, however,
compatible with the elliptic regularity theory;
(ii) We only consider source terms 
in the Lebesgue spaces $L^q(\Dom)$ with $q>2_*:=\frac{2d}{2+d}\ge 1$; 
notice though that this regularity is weaker than
assuming that source terms are in $L^2(\Dom)$, as usually done in the 
literature. 

The traditional approach to tackle the error analysis for
nonconforming approximation techniques are Strang's lemmas. However,
an important shortcoming of this approach whenever the Sobolev smoothness index of the
exact solution is barely larger than one, is that it is not possible
to insert the exact solution in the first argument of the discrete
bilinear form. To do so, one needs to assume some additional
regularity on the exact solution which often goes beyond the
regularity provided by 
the problem at hand.  This approach has nevertheless been used by many
authors to analyze discontinuous Galerkin (dG) methods (see, \eg
\citep{DiPEr:12,ErnGu:06} and the references therein).  One way to
overcome the limitations of Strang's Second Lemma has been proposed by
\citet{Gudi:10}. The key idea consists of introducing a mapping that
transforms the discrete test functions into elements of the exact test
space. An important property of this operator is that its kernel is
composed of discrete (test) functions that are only needed to
``stabilize'' the discrete bilinear form, but do not contribute to the
interpolation properties of the approximation setting. We refer to
this mapping as trimming operator.  The notion of trimming operator
has ben used in \cite{Li_Shipeng_2013} to perform the analysis of the
Crouzeix--Raviart approximation of the diffusion problem 
and source term in $\Ldeux$ (see \eg the
definitions (5)--(7) and the identity (11) therein). The trimmed error
estimate (which is sometimes referred to as ``medius analysis'' in the
literature) has been applied in \cite{Gudi:10} to the IPDG
approximation of the Laplace equation with a source term in $\Ldeux$
and to a fourth-order problem; it has been applied to the Stokes
equations in \cite{BaCGG:14} and to the linear elasticity equations
in~\cite{CarSc:15}.
One problem with methods
using the trimming operator, though, is that they require constructing
$H^1$-conforming discrete quasi-approximation operators that do not
account for the diffusivity contrast; this entails error estimates
with constants that depend on the diffusivity contrast, \ie these
error estimates are not robust.

It is shown in \citep{Ern_Guermond_cmam_2017} in the case of
Nitsche's boundary penalty method that the dependency of
the constants with respect to the diffusivity contrast can be
eliminated by introducing an alternative technique based on
mollification and an extension of the notion of the normal derivative.
The objective of the present paper is to revisit and extend
\citep{Ern_Guermond_cmam_2017}.  The analysis presented here is
significantly simplified and modified to include the Crouzeix--Raviart
approximation, the IPDG approximation, and the HHO
approximation. One key novelty is the introduction of the weighted
bilinear form \eqref{eq:def_n_sharp} that accounts for the
consistency default in all the cases. 
The present analysis hinges on two key
ideas which are now part of the numerical analysis folklore.  To the
best of our knowledge, these ideas have been introduced/used in Lemma~4.7 in
\cite{AmBDG:98}, Lemma~2.3 and Corollary~3.1 in \cite{BerHe:02} and
Lemma~8.2 in \cite{Buffa_Perugia_SINUM_2006}. However, we believe
that detailed proofs are seemingly missing in the literature, and another 
purpose of this paper is to fill this gap.

The first key idea is a
face-to-cell lifting operator. Such an operator is mentioned in
Lemma~4.7 in \citep{AmBDG:98}, and its construction is briefly
discussed.  The weights used in the norms therein, though, cannot give
estimates that are uniform with respect to the mesh size.  This
operator is also mentioned in Lemma~2.3 in \citep{BerHe:02}. The
authors claim that the face-to-cell operator has been constructed in
\cite[Eq.~(5.1)]{Bernardi_Girault_1998}, which is unclear to us. A
similar operator is invoked in Lemma~8.2 in
\citep{Buffa_Perugia_SINUM_2006}. The operator therein
is constructed on the
reference element $\wK$ and its stability properties are proved in the
 Sobolev scale $(H^s(\wK))_{s\in(0,1)}$. The authors invoke also the
Sobolev scale $(H^s(K))_{s\in(0,1)}$ for arbitrary cells $K$ in a mesh
$\calT_h$ belonging the shape-regular sequence $\famTh$. The norm
equipping $H^s(K)$ is not explicitly defined therein, which leads to
one statement that looks questionable (see \eg Eq.~(8.11) therein; a fix has been proposed
in \citep[Lem.~A.3]{Bonito_Guermond_Luddens_M2aN_2016}). In
particular, it is unclear how to keep track of constants that depend on
$K$ when one uses the real interpolation method to define $H^s(K)$.
In order to clarify the status of this face-to-cell operator, which is
essential for our analysis, and without claiming originality, we give
(recall) all the details of its construction in the proof of
Lemma~\ref{lem:lifting_face_to_cell}.  As in
\citep[Lem.~4.7]{AmBDG:98}, we use the Sobolev--Slobodeckij norm to
equip the fractional-order Sobolev spaces; this allows us to track
all the constants easily.  

The second key idea introduced in
the above papers is that of extending the notion of face integrals by
using a duality argument together with the face-to-cell operator. The
argument is deployed in Corollary~3.3 in \citep{BerHe:02}, but the
sketch of the proof has typos (\eg an average has to be removed to
make the inverse estimate in step (1) correct). This corollary is
quoted and invoked in \cite[Lem.~2.1]{CaiHZ:11}; it is the cornerstone
of the argumentation therein.  This argument is also deployed in
Lemma~8.2 in \citep{Buffa_Perugia_SINUM_2006}.  A similar argument is
invoked in \citep[Lem.~4.7]{AmBDG:98} in a slightly different context.
In all the cases one must use a density argument to complete the
proofs, but this argument is omitted and implicitly assumed to hold
true in all the above references. We fill this gap in
Lemma~\ref{lem:identity_n_sharp_dif} and provide the full
argumentation in the proof, including the passage to the limit by
density. The proof invokes mollifiers that commute with differential
operators and behave properly at the boundary of the domain; these
tools have been recently revisited in \citep{ErnGu:16_molli}
elaborating on seminal ideas from \cite{Schoberl_2001}.

\section{Preliminaries}\label{Sec:Preliminaries}
In this section we introduce the model problem and the discrete setting
for the approximation.

\subsection{Model problem}
Let $\Dom$ be a
Lipschitz domain in $\Real^d$, which we assume for simplicity to be a
polyhedron. We consider the following scalar model problem:
\begin{equation} \label{model_elliptic_contrast}
- \div( \lambda \GRAD u ) = f\quad\text{in $\Dom$}, \qquad
\gamma\upg(u) = g\quad\text{on $\front$},
\end{equation}
where $\gamma\upg:\Hun\to H^{\frac12}(\front)$ is the usual trace map
(the superscript $\upg$ refers to the gradient),
and $g\in H^{\frac12}(\front)$ is the Dirichlet boundary data. The
scalar-valued diffusion coefficient $\lambda \in L^\infty(\Dom)$ is
assumed to be uniformly bounded from below away from zero.  For
simplicity, we also assume that $\lambda$ is piecewise constant in
$\Dom$, \ie there is a partition of $\Dom$ into $M$ disjoint Lipschitz
polyhedra $D_1,\cdots,D_M$ \sth $\lambda_{|\Dom_i}$ is a positive real
number for all $i\in\intset{1}{M}$.  

It is standard in the literature to assume that $f\in\Ldeux$. We are going to relax this
hypothesis in this paper by only assuming  that $f\in L^q(\Dom)$ with
$q > \frac{2d}{2+d}$. Note that $q>1$ since $d\ge2$. 
Note also that $L^q(\Dom) \hookrightarrow \Hmun$ since $\Hunz \hookrightarrow H^{q'}(\Dom)$
with the convention that $\frac{1}{q}+\frac{1}{q'}=1$.
Since
$\frac{2d}{2+d}<2$, we are going to assume without loss of generality that 
$q\le 2$.

In the case of the homogeneous Dirichlet condition ($g=0$), 
the weak formulation of the model problem~\eqref{model_elliptic_contrast} 
is as follows: 
\begin{equation}
\left\{
\begin{array}{l}
\text{Find $u \in V:=\Hunz$ such that}\\[2pt]
a(u,w) = \ell(w), \quad \forall w \in V,
\end{array}
\right.
\label{eq:weak_PDE_contrast}
\end{equation}
with the bilinear and linear forms
\begin{equation} \label{eq:a_l_contrast}
a(v,w) := \int_{\Dom} \lambda \GRAD v \SCAL \GRAD w\dif x, \qquad
\ell(w):=\int_\Dom fw\dif x.
\end{equation}
The bilinear form $a$ is coercive in $V$ owing to the
Poincar\'e--Steklov inequality, and it is also bounded on $V\CROSS V$
owing to the Cauchy--Schwarz inequality. The linear form $\ell$ is
bounded on $V$ since the Sobolev embedding theorem and H\"older's
inequality imply that
$|\ell(w)|\le \|f\|_{L^q(\Dom)} \|w\|_{L^{q'}(\Dom)} \le
c\|f\|_{L^q(\Dom)}\|w\|_{H^1(\Dom)}$. Note that $q\ge\frac{2d}{2+d}$ is
the minimal integrability requirement on $f$ for this boundedness
property to hold true. The above coercivity and boundedness properties
combined with the Lax--Milgram Lemma imply
that~\eqref{eq:weak_PDE_contrast} is well-posed.  For the
non-homogeneous Dirichlet boundary condition, one invokes the
surjectivity of the trace map $\gamma\upg$ to infer the existence of a
lifting of $g$, say $u_g\in \Hun$ \sth $\gamma\upg(u_g)=g$, 
and one decomposes the exact
solution as $u=u_g+u_0$ where $u_0\in \Hunz$ solves the weak
problem~\eqref{eq:weak_PDE_contrast} with $\ell(w)$ replaced by
$\ell_g(w)=\ell(w)-a(u_g,w)$.  The weak formulation thus modified is
well-posed since $\ell_g$ is bounded on $\Hunz$.

The notion of diffusive flux, which is defined as follows, will play an important role
in the paper:
\begin{equation}
\bsigma(v):=-\lambda \GRAD v\in\Ldeuxd, \qquad \forall v\in \Hun.
\end{equation} 
We use boldface notation to denote vector-valued functions and 
vectors in $\Real^d$.

\begin{Lem}[Exact solution] \label{lem:p_q_r}
Assume that there exist $r>0$ and $q\in (\frac{2d}{2+d}, 2]$ such that
the exact solution $u$ is in $H^{1+r}(\Dom)$ and the source term $f$ is 
in $L^q(\Dom)$,  then 
\begin{equation} \label{eq:def_V_loS_contrast}
u\in V\loS := \{v\in \Hunz \tq \bsigma(v)\in \bL^p(\Dom), \; \DIV\bsigma(v) \in L^q(\Dom)\},
\end{equation}
for some real number $p>2$.
\end{Lem}

\bproof 
The Sobolev embedding theorem
implies that there is $p>2$ \sth
$\bH^{r}(\Dom)\hookrightarrow \bL^p(\Dom)$.
Indeed, if $2r<d$, we have
$\bH^{r}(\Dom)\hookrightarrow \bL^{s}(\Dom)$ for all
$s\in [2,\frac{2d}{d-2r}]$ and we can take $p=\frac{2d}{d-2r}>2$, 
whereas if $2r\ge d$, we have
$\bH^{r}(\Dom)\hookrightarrow \bH^{\frac{d}{2}}(\Dom) \hookrightarrow
\bL^s(\Dom)$ for all $s\in [2,\infty)$, and we can take any $p>2$.  
The above argument implies that $\GRAD u\in \bL^p(\Dom)$, and since
$\lambda$ is piecewise constant and $\bsigma(u)=-\lambda\GRAD u$, we
have $\bsigma(u)\in \bL^p(\Dom)$.  Moreover, since $\DIV\bsigma(u)=f$
and $f\in L^q(\Dom)$, we have
$\DIV\bsigma(u)\in L^q(\Dom)$. 
\end{proof}

The regularity assumption $u\in H^{1+r}(\Dom)$, $r>0$, is reasonable
owing to the elliptic regularity theory (see
Theorem~3 in \cite{Jochmann_1999}, Lemma~3.2 in
\cite{Bonito_guermond_Luddens_amd2_2013} or \citet{BerVe:00}). In general, one expects that
$r\le\frac12$ whenever $u$ is supported on at least two contiguous subdomains
where $\lambda$ takes different values; otherwise the normal
derivative of $u$ would be continuous across the interface
separating the two subdomains in question, and owing to the discontinuity of
$\lambda$, the normal component of the diffusive flux $\bsigma(u)$ would 
be discontinuous across the interface, which would contradict the
fact that $\bsigma(u)$ has a weak divergence.  It is however
possible that $r>\frac12$ when the exact solution is supported on one
subdomain only. If $r\ge1$, we notice that one necessarily has $f\in \Ldeux$
(since $f_{|\Dom_i} = \lambda_{|\Dom_i}(\Delta u)_{\Dom_i}$ for all
$i\in\intset{1}{M}$), \ie it is legitimate to assume that $q=2$ if $r\ge 1$.
 
\bRem[Extensions] One could also consider lower-order terms
in~\eqref{model_elliptic_contrast}, \eg 
$- \div( \lambda \GRAD u ) + \bbeta\SCAL\GRAD u + \mu u = f$ with
$\bbeta\in\bW^{1,\infty}(\Dom)$ and $\mu\in L^\infty(\Dom)$ \sth
$\mu-\frac12\DIV\bbeta\ge 0$ \ae in $\Dom$ (for simplicity).  The
 error analysis presented in this paper still applies provided
the lower-order terms are not too large, \eg
$\lambda \ge \max(h\|\bbeta\|_{\bL^\infty(\Omega)},
h^2\|\mu\|_{L^\infty(\Dom)})$,
where $h$ denotes the mesh-size. Standard stabilization techniques 
have to be invoked if
the lower-order terms are large when compared to the 
second-order diffusion operator. Furthermore, the error
analysis can be  extended to account  for a
piecewise constant tensor-valued diffusivity $\difT$; then, 
the various constants in the error estimate
depend on the square-root of the anisotropy
ratios measuring the contrast between the largest and the smallest
eigenvalue of $\difT$ in each subdomain $\Dom_i$.  Finally,
one can consider that the diffusion tensor $\difT$ is piecewise smooth
instead of being piecewise constant; a reasonable requirement is that
$\difT_{|\Dom_i}$ is Lipschitz for all $i\in\intset{1}{M}$. This last
extension is, however, less straightforward because
the discrete diffusive flux is no longer a piecewise polynomial function.  
\end{Rem}

\subsection{Discrete setting}
We introduce in this section the discrete setting that we are going to
use to approximate the solution to~\eqref{eq:weak_PDE_contrast}.  Let
$\calT_h$ be a mesh from a shape-regular sequence $\famTh$. Here,
$\calH$ is a countable set with $0$ as unique accumulation point. 
A generic mesh cell is denoted $K\in\calT_h$ and is conventionally
taken to be an open set. 
We also assume that $\calT_h$ covers each of the subdomains
$\{\Dom_i\}_{i\in\intset{1}{M}}$ exactly so that
$\lambda_K:=\lambda_{|K}$ is constant for all $K\in\calT_h$. Let
$\wKPS$ be the reference finite element; we assume that
$\polP_{k,d}\subset \wP \subset W^{k+1,\infty}(\wK)$ for some $k\ge1$.
Here, $\polP_{k,d}$ is the (real) vector space composed of the $d$-variate polynomials 
of degree at most $k$.
For all $K\in\calT_h$, let $\bT_K:\wK\to K$ be the geometric mapping
and let $\mapKg(v)=v\circ \bT_K$ be the pullback by the geometric
mapping.  We introduce the broken finite element space 
\begin{equation}
P\upb_k(\calT_h) = \bset v_h\in
L^\infty(\Dom) \tq v_{h|K}\in P_K,\, \forall K\in\calT_h\eset,
\end{equation}
where $P_K:=(\mapKg)^{-1}(\wP)\subset W^{k+1,\infty}(K)$.  For any
function $v_h\in P\upb_k(\calT_h)$, we define the broken diffusive
flux $\bsigma(v_h)\in \Ldeuxd$ by setting
$\bsigma(v_h)_{|K} := -\lambda_K \GRAD (v_{h|K})$ for all
$K\in\calT_h$.  Upon introducing the notion of broken gradient
$\GRAD_h : W^{1,p}(\calT_h):=\{v\in L^p(\Dom)\st \GRAD (v_{|K}) \in
L^p(K), \ \forall K\in \calT_h\}$ by setting
$(\GRAD_h v)_{|K} :=\GRAD (v_{|K})$ for all $K\in\calT_h$ and all
$v\in W^{1,p}(\calT_h) $, we have $\bsigma(v_h)=-\lambda \GRAD_hv_h$.

For any cell $K\in\calT_h$ we denote by $\bn_K$ the unit normal vector
on $\partial K$ pointing outward.  We denote by $\calFhi$ the
collection of the mesh interfaces and $\calFhb$ the collection of the
mesh faces at the boundary of $\Dom$.  We assume that $\calT_h$ is
oriented in a generation-compatible way, and for each mesh face
$F\in \calFhi\cup \calFhb$, we denote by $\bn_F$ the unit vector
orienting $F$.  For all $F\in\calFhi$, we denote by
$K_l, K_r\in\calT_h$ the two cells \sth
$F=\partial K_l\cap \partial K_r$ and the unit normal vector $\bn_F$
orienting $F$ points from $K_l$ to $K_r$, \ie
$\bn_F=\bn_{K_l}=-\bn_{K_r}$.  For all $F\in\calFh$, let $\calT_F$ be
the collection of the one or two mesh cells sharing $F$. For all
$K\in\calT_h$, let $\calFK$ be the collection of the faces of $K$ and
let $\epsilon_{K,F}=\bn_F\SCAL\bn_K=\pm1$.  The jump across $F\in\calFhi$ of any
function $v\in W^{1,1}(\calT_h)$ is defined by setting
$\jump{v}_F(\bx) = v_{|K_l}(\bx)-v_{|K_r}(\bx)$ for \ae $\bx\in F$.
If $F\in \calFhb$, this jump is conventionally defined as the trace
on $F$, \ie $\jump{v}_F(\bx) = v_{|K_l}(\bx)$ where 
$F=\partial K_l\cap\partial\Dom$. We omit the subscript $_F$ 
in the jump whenever the context is unambiguous.

\section{The bilinear form $n_\sharp$} \label{Sec:nsharp}

In this section, we give a proper meaning to the normal trace of the
diffusive flux of the solution to \eqref{eq:weak_PDE_contrast} over
each mesh face. The material presented in
\S\ref{sec:face-to-cell_extension} and \S\ref{sec:extension_contrast}
has been introduced in \citep[\S5.3]{Ern_Guermond_cmam_2017} and is
inspired from \cite[Lem.~4.7]{AmBDG:98}, \cite[Cor 3.3]{BerHe:02}, and 
\cite[Lem.~8.2]{Buffa_Perugia_SINUM_2006};
it is included here for the sake of completeness.
The reader familiar with these techniques is invited to jump to
\S\ref{Section:nsharp} where the weighted bilinear form $n_\sharp$ is
introduced.  This bilinear form is the main tool for the
error analysis in \S\ref{Sec:applications}.

\subsection{Face-to-cell lifting operator}
\label{sec:face-to-cell_extension}

Let us first motivate our approach informally. 
Let $K\in\calT_h$ be a mesh cell, let $\calF_K$ be the 
collection of all the faces of $K$,
and let $F\in\calF_K$ be a face of $K$.
Let $\bv$ be a vector field defined on $K$. We are looking for (mild)
regularity requirements on the field $\bv$ 
to give a meaning to the quantity 
$\int_F (\bv\SCAL\bn_K)\phi\dif s$, where $\phi$ is a given smooth function
on $F$ (\eg a polynomial function). 
It is well established that it is possible
to give a weak meaning in $H^{-\frac12}(\partial K)$ to the 
normal trace of $\bv$ on $\partial K$ by means of an
integration by parts formula if $\bv \in \bH(\text{\rm div};K)
:= \{\bv\in \bL^2(K) \tq \DIV\bv \in L^2(K)\}$. In this situation, 
one can define the normal
trace $\gamma_{\partial K}\upd(\bv)\in H^{-\frac12}(\partial K)$ by setting
\begin{equation}
\langle \gamma_{\partial K}\upd(\bv),\psi\rangle_{\partial K} := \int_K \Big( 
\bv\SCAL \GRAD w(\psi) + (\DIV\bv)w(\psi)\Big)\dif x,
\end{equation} 
for all $\psi\in H^{\frac12}(\partial K)$, where $w(\psi)\in H^1(K)$ 
is a lifting of $\psi$, \ie $\gamma\upg_{\partial K}(w(\psi))
= \psi$, and $\gamma\upg_{\partial K}:H^1(K)\to H^{\frac12}(\partial K)$ 
is the trace map locally in $K$. 
Then, one has $\gamma_{\partial K}\upd(\bv)
= \bv_{|\partial K}\SCAL\bn_K$ whenever
$\bv$ is smooth, \eg if $\bv\in \bH(\text{\rm div};K)\cap \bC^0(\overline K)$.
However, the above meaning
is too weak for our purpose because we need to localize the action of the
normal trace to functions $\phi$ only defined on a face $F$, \ie $\phi$ may not be defined over the
whole boundary $\partial K$. The key to achieve this is to extend
$\phi$ by zero from $F$ to $\partial K$. This obliges us to change
the functional setting since the extended function is
no longer in $H^{\frac12}(\partial K)$. In what follows,
we are going to use the fact that the 
zero-extension of a smooth function defined on a face $F$
of $\partial K$ is in $W^{1-\frac1t,t}(\partial K)$ if $t<2$, \ie
$t(1-\frac1t)<1$. Let us now present 
a rigorous construction.
 
Let $p,q$ be two real numbers such that
\begin{equation}
p>2, \qquad q>\frac{2d}{2+d}.\label{eq:p_q_face_to_cell}
\end{equation}
Notice that $q>1$ since $d\ge2$.  Let $\tp\in (2,p]$ be
such that $q\ge \frac{\tp d}{\tp +d}$; this is indeed possible since
$p>2$, $q>\frac{2d}{2+d}$, and the function $z\mapsto \frac{zd}{z+d}$ is increasing over
$\Real_+$. Lemma~\ref{lem:lifting_face_to_cell} shows that there exists
a bounded lifting operator 
\begin{equation}
L_F^K : W^{\frac{1}{\tp},\tp'}(F) \longrightarrow W^{1,\tp'}(K),
\end{equation}
with conjugate number $\tp'$ s.t.~$\frac{1}{\tp}+\frac{1}{\tp'}=1$, so
that for any $\phi\in W^{\frac{1}{\tp},\tp'}(F)$, $L_F^K(\phi)$ is a 
lifting of the zero-extension of $\phi$ to $\partial K$, \ie
\begin{equation}
\label{eq:traces_lifting_p_q}
\gamma\upg_{\partial K}(L_F^K(\phi))_{|\partial K{\setminus} F}=0, \qquad 
\gamma\upg_{\partial K}(L_F^K(\phi))_{|F} =\phi.
\end{equation}
Notice that the domain of $L_F^K$ is of the form $W^{1-\frac1t,t}(F)$ with $t=\tp'<2$, 
which is consistent with the above observation regarding the zero-extension to $\partial K$ 
of functions defined on $F$. We also observe that
\begin{equation}
L_F^K(\phi) \in W^{1,p'}(K) \cap L^{q'}(K),
\end{equation} 
with conjugate numbers $p',q'$ s.t.~$\frac{1}{p}+\frac{1}{p'}=1$,
$\frac{1}{q}+\frac{1}{q'}=1$. Indeed, $L_F^K(\phi) \in W^{1,p'}(K)$
just follows from $p'\le \tp'$ (\ie $\tp\le p$), whereas
$L_F^K(\phi) \in L^{q'}(K)$ follows from 
$W^{1,\tp'}(K) \hookrightarrow L^{q'}(K)$ owing to
the Sobolev Embedding Theorem
(since $q'\le \frac{\tp'd}{d-\tp'}$, as can be verified from
$d\ge 2 > \tp'$ and 
$\frac{1}{\tp'}-\frac{1}{d} = 1-(\frac1\tp+\frac1d)\le
1-\frac1q=\frac{1}{q'}$ because $q\ge\frac{\tp d}{\tp+d}$).
We now state our main result on the 
lifting operator $L_F^K$.

\begin{Lem}[Face-to-cell lifting] \label{lem:lifting_face_to_cell}
Let $p$ and $q$ satisfy \eqref{eq:p_q_face_to_cell}.
Let $\tp\in (2,p]$ be such that $q\ge \frac{\tp d}{\tp +d}$.  
Let $K\in\calT_h$ be a mesh cell and let $F\in\calFK$ be a face of $K$.
There exists a lifting operator
$L_F^K : W^{\frac{1}{\tp},\tp'}(F) \rightarrow W^{1,\tp'}(K)$
satisfying~\eqref{eq:traces_lifting_p_q}, and there exists
$c$, uniform \wrt $\hinH$, but depending on $p$ and $q$,
\sth the following holds true:
\begin{equation} 
\label{eq:stab_lifting_p_q}
h_K^{\frac{d}{p}}|L_F^K(\phi)|_{W^{1,p'}(K)} 
+ h_K^{-1+\frac{d}{q}}\|L_F^K(\phi)\|_{L^{q'}(K)}
\le c\, h_K^{-\frac{1}{\tp} + \frac{d}{\tp}}\|\phi\|_{W^{\frac{1}{\tp},\tp'}(F)}, 
\end{equation}
for all $\phi\in W^{\frac{1}{\tp},\tp'}(F)$ 
with the norm $\|\phi\|_{W^{\frac{1}{\tp},\tp'}(F)} := \|\phi\|_{L^{\tp'}(F)}+ h_F^{\frac{1}{\tp}}
|\phi|_{W^{\frac{1}{\tp},\tp'}(F)}$.
\end{Lem}

\bproof
(1) The face-to-cell lifting operator $L_F^K$ is constructed from a
lifting operator $L_{\wF}^{\wK}$ on the reference cell.
Let $\wK$ be the reference cell and let $\wF$ be one of its faces.  
Let us define the operator
$L_{\wF}^{\wK} : W^{\frac{1}{\tp},\tp'}(\wF)\rightarrow W^{1,\tp'}(\wK)$. 
For any function
$\psi \in W^{\frac{1}{\tp},\tp'}(\wF)$, 
let $\widetilde\psi$ denote the zero-extension
of $\psi$ to $\partial \wK$.  
Owing to \citet[Thm.~1.4.2.4, Cor.~1.4.4.5]{Gr85}, $\widetilde \psi$ is in
$W^{\frac{1}{\tp},\tp'}(\partial\wK)$ 
since $\frac{\tp'}{\tp} = \frac{1}{\tp-1} <1$ (\ie $\tp>2$), and we have
$\|\widetilde \psi\|_{W^{\frac{1}{\tp},\tp'}(\partial\wK)}\le \wc_1
\|\psi\|_{W^{\frac{1}{\tp},\tp'}(\wF)}$ with the norm 
$\|\psi\|_{W^{\frac{1}{\tp},\tp'}(\wF)} := 
\|\psi\|_{L^{\tp'}(\wF)}+ \ell_\wK^{\frac{1}{\tp}}
|\psi|_{W^{\frac{1}{\tp},\tp'}(\wF)}$ where $\ell_\wK=1$ is a 
length scale associated with $\wK$. 
Then we use the surjectivity of the trace map
$\gamma\upg_{\wK}: W^{1,\tp'}(\wK)\to W^{\frac{1}{\tp},\tp'}(\partial\wK)$
(see \citet[Thm.~1.I]{Gagliardo:57}) to define
$L_{\wF}^{\wK}(\psi) \in W^{1,\tp'}(\wK)$ \sth
$\gamma\upg_{\wK}(L_{\wF}^{\wK}(\psi)) = \widetilde\psi$ and
$\|L_{\wF}^{\wK}(\psi)\|_{W^{1,\tp'}(\wK)} \le \wc_2
\|\widetilde\psi\|_{W^{\frac{1}{\tp},\tp'}(\partial\wK)}$, \ie
$\|L_{\wF}^{\wK}(\psi)\|_{W^{1,\tp'}(\wK)} \le \wc
\|\psi\|_{W^{\frac{1}{\tp},\tp'}(\wF)}$,
with $\wc=\wc_1\wc_2$. By construction, we have 
$\gamma\upg_{\partial \wK}(L_{\wF}^{\wK}(\psi))_{|\wF}=\psi$ and 
$\gamma\upg_{\partial \wK}(L_{\wF}^{\wK}(\psi))_{|\partial\wK\setminus\wF}=0$.
\\
(2) We define the lifting operator 
$L_F^K:W^{\frac{1}{\tp},\tp'}(F) \rightarrow W^{1,\tp'}(K)$ by setting
\begin{equation}
L_F^K(\phi)(\bx) := 
L_{\wF}^{\wK}( \phi\circ \bT_{K|\wF}) (\bT_K^{-1}(\bx)), \quad \forall \bx\in K,\quad
 \forall \phi\in W^{\frac{1}{\tp},\tp'}(F), \label{def_of_LFK_from_LwFWK}
\end{equation}
where $\bT_K :\wK\to K$ is the geometric mapping and
$\wF=\bT_K^{-1}(F)$. By definition, if $\bx \in F$, 
then $\wbx :=\bT_K^{-1}(\bx)\in \wF$ and
$\bT_{K|\wF}(\wbx)=\bx$, so that 
\[
\gamma\upg_{\partial K}(L_F^K(\phi))(\bx) =
\gamma\upg_{\partial \wK}(L_{\wF}^{\wK}(\phi\circ \bT_{K|\wF}))(\wbx)=\phi(\bT_{K|\wF}(\wbx))=\phi(\bx),
\]
whereas if $\bx\in\partial K\setminus F$, then
$\wbx\in \partial\wK\setminus \wF$, so that
$\gamma\upg_{\partial \wK}(L_{\wF}^{\wK}(\phi\circ
\bT_{K|\wF}))(\wbx)=0$. The above argument shows
that~\eqref{eq:traces_lifting_p_q} holds true.
\\
(3) It remains to prove~\eqref{eq:stab_lifting_p_q}. Let us first
bound $|L_F^K(\phi)|_{W^{1,p'}(K)}$. Notice that the definition
of $L_K^F$ is equivalent to 
$L_F^K(\phi)\circ \bT_K(\wbx) := L_{\wF}^{\wK}( \phi\circ
\bT_{K|\wF}) (\wbx)$; that is,
$\mapKg(L_F^K(\phi)) := L_{\wF}^{\wK}(\mapFg(\phi))$, where $\mapKg$
is the pullback by $\bT_K$, and $\mapFg$ is the pullback by
$\bT_{K|\wF}$.  Denoting by $\Jac_K$ the Jacobian of the geometric mapping $\bT_K$, we infer that
\begin{align*}
|L_F^K(\phi)|_{W^{1,p'}(K)} &\le c\, \|\Jac_K^{-1}\|_{\ell^2} |\det(\Jac_K)|^{\frac{1}{p'}} |L_{\wF}^{\wK}(\mapFg(\phi))|_{W^{1,p'}(\wK)} \\
&\le c'\, \|\Jac_K^{-1}\|_{\ell^2} |\det(\Jac_K)|^{\frac{1}{p'}} |L_{\wF}^{\wK}(\mapFg(\phi))|_{W^{1,\tp'}(\wK)} \\
&\le c''\, \|\Jac_K^{-1}\|_{\ell^2} |\det(\Jac_K)|^{\frac{1}{p'}} \|\mapFg(\phi)\|_{W^{\frac{1}{\tp},\tp'}(\wF)},
\end{align*}
where the first inequality follows from the chain rule, 
the second is a consequence of $\tp'\ge p'$ (since $\tp\le p$), and
the third follows from the stability of the reference lifting operator
$L_{\wF}^{\wK}$. Using now the chain rule and the shape-regularity of the mesh sequence, 
we infer that
$\|\mapFg(\phi)\|_{W^{\frac{1}{\tp},\tp'}(\wF)} \le c|\det(\Jac_F)|^{-\frac{1}{\tp'}}
\|\phi\|_{W^{\frac{1}{\tp},\tp'}(F)}$,
where $\Jac_F$ is the Jacobian of the mapping $\bT_{K|\wF} : \wF\to
F$. Combining these bounds, we obtain
\begin{align*} 
|L_F^K(\phi)|_{W^{1,p'}(K)} &\le c\, \|\Jac_K^{-1}\|_{\ell^2} |\det(\Jac_K)|^{\frac{1}{p'}} |\det(\Jac_F)|^{-\frac{1}{\tp'}} \|\phi\|_{W^{\frac{1}{\tp},\tp'}(F)} \\
&\le c'\, h_K^{-\frac{1}{\tp}+d(\frac{1}{\tp}-\frac1p)} \|\phi\|_{W^{\frac{1}{\tp},\tp'}(F)},
\end{align*}
where the second bound follows from the shape-regularity of the mesh sequence. 
This proves the bound on $|L_F^K(\phi)|_{W^{1,p'}(K)}$ 
in~\eqref{eq:stab_lifting_p_q}. The proof of the bound on 
$\|L_F^K(\phi)\|_{L^{q'}(K)}$ uses similar arguments 
together with $W^{1,\tp'}(\wK) \hookrightarrow L^{q'}(\wK)$
owing to the Sobolev Embedding Theorem
and $q'\le \frac{\tp'd}{d-\tp'}$ (as already shown above).
\end{proof}

\subsection{Face localization of the normal diffusive flux}
\label{sec:extension_contrast}

Let $K\in\calT_h$ be a mesh cell, $F\in\calFK$ be a face of $K$, and
consider the following functional space:
\begin{equation} 
\bS\upd(K):=\bset \btau\in \bL^p(K)\tq \div \btau \in L^q(K)\eset,
\label{def_of_Vd_contrasted_diffusivity}
\end{equation}
equipped with the 
following dimensionally-consistent norm:
\begin{equation}
\|\btau\|_{\bS\upd(K)}:=\|\btau\|_{\bL^p(K)} +
  h_{K}^{1+d(\frac1p-\frac1q)}\|\DIV\btau\|_{L^q(K)}. \label{def_of_bVupd_K}
\end{equation}
With the lifting operator $L_F^K$ in hand, we now define the normal trace on the face $F$ of $K$
of any field $\btau\in \bS\upd(K)$ to be the linear form in $(W^{\frac{1}{\tp},\tp'}(F))'$
denoted by $(\btau\SCAL\bn_K)_{|F}$ and whose action on any function $\phi \in W^{\frac{1}{\tp},\tp'}(F)$ is 
\begin{equation} \label{eq:def_comp_normale} 
\langle (\btau\SCAL\bn_K)_{|F},\phi \rangle_{F} := \int_K \Big(\btau\SCAL \GRAD
L_F^K(\phi) + (\DIV\btau) L_F^K(\phi)\Big) \dif x.
\end{equation}
Here, $\langle\SCAL,\SCAL\rangle_F$ denotes
the duality pairing between $(W^{\frac{1}{\tp},\tp'}(F))'$ 
and $W^{\frac{1}{\tp},\tp'}(F)$. Notice that
the right-hand side of~\eqref{eq:def_comp_normale} is
well-defined owing to H\"older's inequality
and~\eqref{eq:stab_lifting_p_q}. 
Owing to~\eqref{eq:traces_lifting_p_q}, we readily verify 
that we have indeed defined an extension of the normal trace 
since we have
$\langle (\btau\SCAL\bn_K)_{|F},\phi \rangle_F= 
\int_F (\btau\SCAL\bn_K)\phi \dif s$ whenever the field $\btau$ is smooth.
Let us now derive an important bound on the linear form
$(\btau\SCAL\bn_K)_{|F}$ when acting on a function from the space
$P_F$, which we define to be composed of the restrictions to $F$ 
of the functions in $P_K$. 
Note that $P_F \subset W^{\frac{1}{\tp},\tp'}(F)$.

\begin{Lem}[Bound on normal component]\label{lem:bnd_normal_contrast}
There exists a constant $c$, uniform \wrt $\hinH$,
but depending on $p$ and $q$, 
\sth the following holds true: 
\begin{align}
  \label{eq:bnd_normal_dual_h} 
  |\langle(\btau\SCAL\bn_K)_{|F},\phi_h \rangle_F| &\le c\,h_{K}^{d(\frac12-\frac1p)}\|\btau\|_{\bS\upd(K)}  h_F^{-\frac12}\|\phi_h\|_{L^2(F)},
\end{align}%
for all $\btau\in \bS\upd(K)$, all $\phi_h\in P_F$, 
all $K\in\calT_h$, and all $F\in\calF_K$.
\end{Lem}

\bproof
A direct consequence of~\eqref{eq:def_comp_normale}, 
H\"older's inequality, and Lemma~\ref{lem:lifting_face_to_cell} is that
\[
|\langle (\btau\SCAL\bn_K)_{|F},\phi\rangle_F|\le c\,h_K^{-\frac{1}{\tp}+d(\frac{1}{\tp}-\frac{1}{p})}  \|\btau\|_{\bS\upd(K)} \|\phi\|_{W^{\frac{1}{\tp},\tp'}(F)},
\]
for all $\phi\in W^{\frac{1}{\tp},\tp'}(F)$. Recalling that 
$\|\phi\|_{W^{\frac{1}{\tp},\tp'}(F)} = \|\phi\|_{L^{\tp'}(F)}+ h_F^{\frac{1}{\tp}}
|\phi|_{W^{\frac{1}{\tp},\tp'}(F)}$, the shape-regularity of the mesh sequence implies that the following inverse inequality 
$\|\phi_h\|_{W^{\frac{1}{\tp},\tp'}(F)} \le c
h_F^{(d-1)(\frac12-\frac{1}{\tp})}\|\phi_h\|_{L^2(F)}$ holds true for all $\phi_h\in P_F$
(note that $\frac12-\frac{1}{\tp} = \frac{1}{\tp'}-\frac12$).
The estimate~\eqref{eq:bnd_normal_dual_h} follows readily.
\end{proof}

\subsection{Definition of $n_\sharp$ and key identities} \label{Section:nsharp}

Let us consider the functional space $V\loS$ 
defined in~\eqref{eq:def_V_loS_contrast}. For all $v\in V\loS$, 
Lemma~\ref{lem:p_q_r} shows that
$\bsigma(v)_{|K} \in \bS\upd(K)$ for all $K\in\calT_h$,
and Lemma~\ref{lem:bnd_normal_contrast} 
implies that  it is
possible to give a meaning by duality to the normal component of
$\bsigma(v)_{|K}$ on all the faces of $K$ separately.
Moreover, since we have set $\bsigma(v_h)_{|K}:=-\lambda_K
\GRAD (v_{h|K})$ for all $v_h\in P\upb_k(\calT_h)$, and since we have
$P_K\subset W^{k+1,\infty}(K)$ with $k\ge1$, we infer that $\bsigma(v_h)_{|K}
\in \bS\upd(K)$ as well. Thus, $\bsigma(v)_{|K}\in \bS\upd(K)$ for all 
$v\in (V\loS+P\upb_k(\calT_h))$.
Let us now introduce the bilinear form 
$n_\sharp:(V\loS+P\upb_k(\calT_h)) \times P\upb_k(\calT_h)\to \Real$ defined as
follows:
\begin{align} \label{eq:def_n_sharp}
n_\sharp(v,w_h) := {}&\sum_{F\in\calFh} \sum_{K\in\calT_F}
\epsilon_{K,F}\theta_{K,F} 
\langle (\bsigma(v)_{|K}\SCAL\bn_K)_{|F},\jump{w_h} \rangle_F,
\end{align}
where the weights $\theta_{K,F}$ are still unspecified but are
assumed to satisfy
\begin{equation}
\theta_{K_l,F},\theta_{K_r,F}\in [0,1]\quad \text{and} \quad
\theta_{K_l,F}+\theta_{K_r,F}=1,
\qquad \forall F\in\calFhi, \label{convex_interface_weights}
\end{equation} 
whereas for all $F\in\calFhb$ with
$F=\partial K_l\cap \front$, we set
$\theta_{K_l,F}:=1$, $\theta_{K_r,F}=:0$.
We will see in~\eqref{eq:def_weight_theta} below 
how these weights must depend on the
diffusion coefficient to get a robust boundedness estimate on $n_\sharp$.
The definition~\eqref{eq:def_n_sharp} is meaningful since
$\jump{w_h}_F \in P_F$ for all $w_h\in P\upb_k(\calT_h)$. 
The factor $\epsilon_{K,F}$ in~\eqref{eq:def_n_sharp}
handles the relative orientation 
of $\bn_K$ and $\bn_F$.
For all $v\in W^{1,1}(\calT_h)$, we
define weighted averages as follows for \ae $\bx\in F\in\calFhi$:
\begin{subequations} \label{eq:def_weight_avg} \begin{align}
\avg{v}_{F,\theta}(\bx)&:= \theta_{K_l,F} v_{|K_l}(\bx) + \theta_{K_r,F} v_{|K_r}(\bx), \\
\avg{v}_{F,\bar\theta}(\bx)&:= \theta_{K_r,F} v_{|K_l}(\bx) + \theta_{K_l,F} v_{|K_r}(\bx).
\end{align}\end{subequations}
Whenever $\theta_{K_l,F}=\theta_{K_r,F}=\frac12$, these two definitions
coincide with the usual arithmetic average. 
On boundary faces $F\in\calFhb$, we have
$\avg{v}_{F,\theta}(\bx)=v_{|K_l}(\bx)$, and $\avg{v}_{F,\bar\theta}(\bx)=0$ 
for \ae $\bx\in F$.
We omit the subscript $_F$ 
whenever the context is unambiguous.  
The following
identity will be useful: 
\begin{equation} \label{eq:identity_jump_wei_avg}
\jump{vw} = \avg{v}_{\theta} \jump{w} + \jump{v}\avg{w}_{\bar\theta}.
\end{equation}

The following lemma is fundamental to understand the role that the
bilinear form $n_\sharp$ will play in the next section in the
analysis of various nonconforming approximation methods.
\begin{Lem}[Identities for $n_\sharp$] \label{lem:identity_n_sharp_dif}
The following holds true for any choice of weights $\{\theta_{K,F}\}_{F\in \calFh,K\in\calT_F}$ 
and for all $w_h\in P\upb_k(\calT_h)$, all $v_h\in P\upb_k(\calT_h)$, and all $v\in V\loS$:
\begin{subequations}\begin{align} \label{eq:n_sharp_1}
n_\sharp(v_h,w_h) &= \sum_{F\in\calFh} \int_F \avg{\bsigma(v_h)}_\theta\SCAL\bn_F \jump{w_h} \dif s,\\
\label{eq:n_sharp_2}
n_\sharp(v,w_h) &= \sum_{K\in\calT_h} \int_K \Big(\bsigma(v)\SCAL\GRAD w_{h|K}+
 (\DIV\bsigma(v)) w_{h|K}\Big) \dif x.
\end{align}\end{subequations}
\end{Lem}

\bproof (1) Proof of~\eqref{eq:n_sharp_1}. Let
$v_h,w_h\in P\upb_k(\calT_h)$. Since the restriction of $\bsigma(v_h)$
to each mesh cell is smooth, and since the restriction of
$L_F^K(\jump{w_h})$ to $\partial K$ is nonzero only on the face $F\in\calFK$ 
where it coincides with $\jump{w_h}$, we have
\begin{align*}
\langle (\bsigma(v_h)_{|K}\SCAL\bn_K)_{|F},\jump{w_h} \rangle_F 
&= \int_K \Big( \bsigma(v_h)_{|K}\SCAL \GRAD L_F^K(\jump{w_h}) 
+ (\DIV\bsigma(v_h)_{|K})L_F^K(\jump{w_h})\Big) \dif x \\ 
&=\int_{\partial K} \bsigma(v_h)_{|K}\SCAL\bn_KL_F^K(\jump{w_h})\dif s
= \int_F \bsigma(v_h)_{|K}\SCAL\bn_K\jump{w_h}\dif s,
\end{align*}
where we used the divergence formula in $K$. Therefore,
after using the definitions of $\epsilon_{K,F}$ and of $\theta_{K,F}$, 
we obtain 
\begin{align*}
n_\sharp(v_h,w_h) &= \sum_{F\in\calFh} \sum_{K\in\calT_F}
\epsilon_{K,F} \theta_{K,F} \int_F \bsigma(v_h)_{|K}\SCAL\bn_K\jump{w_h}\dif s \\
&= \sum_{F\in\calFh} \int_F \avg{\bsigma(v_h)}_\theta\SCAL\bn_F \jump{w_h} \dif s.
\end{align*} 
(2) Proof of~\eqref{eq:n_sharp_2}. Let $v\in V\loS$ and $w_h\in P\upb_k(\calT_h)$.
Let $\calK_\delta\upd:\bL^1(\Dom) \to \bC^\infty(\overline\Dom)$ and
$\calK_\delta\upb:L^1(\Dom) \to C^\infty(\overline\Dom)$ be the mollification
operators introduced in~\citep[\S3.2]{ErnGu:16_molli}. 
These two operators satisfy the following key commuting property:
\begin{equation} \label{eq:commut_contrast}
\DIV(\calK_\delta\upd(\btau)) = \calK_\delta\upb(\DIV\btau),
\end{equation}
for all $\btau\in \bL^1(\Dom)$ s.t.~$\DIV\btau\in L^1(\Dom)$.  It is
important to realize that this property can be applied to $\bsigma(v)$
for all $v\in V\loS$ since $\DIV\bsigma(v)\in L^1(\Dom)$ by definition
of $V\loS$. (Note that this property cannot be applied to
$\bsigma(v_h)$ with $v_h\in P\upb_k(\calT_h)$, since the normal component of
$\bsigma(v_h)$ is in general discontinuous across
the mesh interfaces, \ie $\bsigma(v_h)$ does not have a 
weak divergence.) Let us consider the mollified
bilinear form
\[
n_{\sharp\delta}(v,w_h) := \sum_{F\in\calFh} \sum_{K\in\calT_F}
\epsilon_{K,F}\theta_{K,F} \langle (\calK_\delta\upd(\bsigma(v))_{|K}\SCAL\bn_K)_{|F},\jump{w_h} \rangle_F.
\]
Owing to the commuting property~\eqref{eq:commut_contrast}, we infer that
\begin{multline*}
\langle (\calK_\delta\upd(\bsigma(v))_{|K}\SCAL\bn_K)_{|F},\jump{w_h} \rangle_F = \\ \int_K \Big(  \calK_\delta\upd(\bsigma(v))\SCAL L_F^K(\jump{w_h}) 
+ \calK_\delta\upb(\DIV\bsigma(v))L_F^K(\jump{w_h}) \Big) \dif x.
\end{multline*}
Then Theorem~3.3 from \citep{ErnGu:16_molli} implies that
\begin{multline*}
\lim_{\delta\to 0}\int_K \Big(  \calK_\delta\upd(\bsigma(v))\SCAL L_F^K(\jump{w_h}) 
+ \calK_\delta\upb((\DIV\bsigma(v)))L_F^K(\jump{w_h}) \Big) \dif x = \\
\int_K \Big( \bsigma(v)\SCAL L_F^K(\jump{w_h}) 
+ (\DIV\bsigma(v))L_F^K(\jump{w_h}) \Big) \dif x 
= \langle (\bsigma(v)_{|K}\SCAL\bn_K)_{|F},\jump{w_h} \rangle_F.
\end{multline*}
Summing over the mesh faces and the associated mesh cells, we infer that
\[
\lim_{\delta\to 0} n_{\sharp\delta}(v,w_h) = n_{\sharp}(v,w_h).
\]
Moreover, since the mollified function $\calK_\delta\upd(\bsigma(v))$
is smooth, by repeating the calculation done in Step~(1), we also have
\begin{align*}
n_{\sharp\delta}(v,w_h) 
&= \sum_{F\in\calFh} \int_F \avg{\calK_\delta\upd(\bsigma(v))}_\theta\SCAL\bn_F \jump{w_h}\dif s.
\end{align*}
Using the identity~\eqref{eq:identity_jump_wei_avg} with
$\jump{\calK_\delta\upd(\bsigma(v))}\SCAL\bn_F = 0$ for all
$F\in\calFhi$, recalling that
$\jump{w_h \calK_\delta\upd(\bsigma(v))}=w_h
\calK_\delta\upd(\bsigma(v))_{|F}$ for all $F\in\calFhb$, and using the
divergence formula in $K$ and the commuting
property~\eqref{eq:commut_contrast}, we obtain
\begin{align*}
n_{\sharp\delta}(v,w_h) &= \sum_{F\in\calFh} \int_F \avg{\calK_\delta\upd(\bsigma(v))}_\theta\SCAL\bn_F \jump{w_h}\dif s + \sum_{F\in\calFhi} \int_F \jump{\calK_\delta\upd(\bsigma(v))}\SCAL\bn_F\avg{w_h}_{\bar\theta}\dif s \\
& = \sum_{F\in\calFh} \int_F \jump{w_h \calK_\delta\upd(\bsigma(v))}\SCAL\bn_F\dif s = 
\sum_{K\in\calT_h} \int_{\partial K} \calK_\delta\upd(\bsigma(v))\SCAL\bn_K w_{h|K}\dif s \\
&= \sum_{K\in\calT_h} \int_K \Big( \calK_\delta\upd(\bsigma(v))\SCAL\GRAD w_{h|K} +
 \calK_\delta\upb(\DIV\bsigma(v))w_{h|K}\Big) \dif x.
\end{align*}
Invoking again Theorem~3.3 from \citep{ErnGu:16_molli} leads to the assertion since
\[
\lim_{\delta\to 0} n_{\sharp\delta}(v,w_h) = \sum_{K\in\calT_h} \int_K \Big(\bsigma(v)\SCAL\GRAD w_{h|K}+
 (\DIV\bsigma(v))w_{h|K}\Big) \dif x. 
\]
\end{proof}

\bRem[Identity \eqref{eq:n_sharp_2}] 
The identity~\eqref{eq:n_sharp_2} is the key tool 
to assert in a weak sense that $\bsigma(v)\SCAL \bn$ is continuous
across the mesh interfaces without the need to assume 
that $v$ is smooth, say $v\in H^{1+r}(\Dom)$ with $r>\frac12$.
\end{Rem}

We now establish an important boundedness estimate on the bilinear
form $n_\sharp$. Since $\bsigma(v)_{|K}\in \bS\upd(K)$ for all $K\in\calT_h$ 
and all $v\in V\loS+P\upb_k(\calT_h)$, we can equip the space 
$V\loS+P\upb_k(\calT_h)$ with the seminorm
\begin{equation}
|v|_{n_\sharp}^2 
:= \sum_{K\in \calT_h} \lambda_K^{-1}\Big(h_K^{2d(\frac12-\frac1p)} \|\bsigma(v)_{|K}\|_{\bL^p(K)}^2
+h_K^{2d(\frac{2+d}{2d}-\frac{1}{q})} 
\|\DIV\bsigma(v)_{|K}\|_{L^q(K)}^2\Big).\label{eq:seminorm_sharp_n}
\end{equation}%
We notice that this seminorm is dimensionally-consistent with the classical
energy-norm defined as $\sum_{K\in\calT_h} \lambda_K\|\GRAD v_{|K}\|_{\bL^2(K)}^2$.
Straightforward algebra shows that $|v|_\sharp \le c\lambda_\flat^{-\frac12}
( \ell_\Dom^{d(\frac12-\frac1p)} 
\|\bsigma(v)\|_{\bL^p(\Dom)}
+ \ell_\Dom^{d(\frac{2+d}{2d}-\frac{1}{q})} \|\DIV\bsigma(v)\|_{L^q(\Dom)})$,
for all $v\in V\loS$; here $\ell_\Dom$ denotes a characteristic length of $\Dom$.
(Recall that $\|a\|_{\ell^s(\calI)}\le \|a\|_{\ell^t(\calI)}$ for any finite sequence 
$(a_i)_{i\in\calI}$ if $0<t\le s$, and we assumed that $q\le 2$.)

In order to get robust error estimates with respect to $\lambda$, it
is important to avoid any 
dependency on the ratio of the values taken by $\lambda$ in two adjacent
subdomains; otherwise, the error estimates become meaningless when the
diffusion coefficient $\lambda$ is highly contrasted.  To avoid such
dependencies, we introduce the following
diffusion-dependent
weights for all $F\in\calFhi$, with $F=\partial K_l\cap \partial K_r$:
\begin{equation} \label{eq:def_weight_theta}
\theta_{K_l,F} := \frac{\lambda_{K_r}}{\lambda_{K_l}+\lambda_{K_r}}, \qquad
\theta_{K_r,F} := \frac{\lambda_{K_l}}{\lambda_{K_l}+\lambda_{K_r}}.
\end{equation} 
We also define 
\begin{equation}
\lambda_F := \frac{2\lambda_{K_l}\lambda_{K_r}}{\lambda_{K_l}+\lambda_{K_r}}\ \text{if $F\in\calFhi$} 
\quad\text{and}\quad \lambda_F:=\lambda_{K_{l}} \ \text{if $F\in\calFhb$}.
\label{def_of_lambdaF}
\end{equation}
The two properties we are going to use are
that $\mes{\calT_F}\lambda_{K}\theta_{K,F}=\lambda_{F}$, 
for all $K\in\calT_F$, and $\lambda_F\le \min_{K\in \calT_F} \lambda_K$.
(Here $\mes{\calT_F}$ denotes the cardinality of $\calT_F$.)

\begin{Lem}[Boundedness of $n_\sharp$] \label{lem:bnd_n_sharp} 
With the weights defined in \eqref{eq:def_weight_theta} and 
$\lambda_F$ defined in \eqref{def_of_lambdaF} for all $F\in\calFh$,
there is $c$, uniform \wrt $\hinH$ and $\lambda$, but depending on
$p$ and $q$, \sth
the following holds true for all
$v\in V\loS+P\upb_k(\calT_h)$ and all $w_h\in P\upb_k(\calT_h)$:
\begin{equation}
|n_\sharp(v,w_h)| \le c\, |v|_{n_\sharp} \bigg( \sum_{F\in\calFh} 
\lambda_F h_F^{-1} \|\jump{w_h}\|_{L^2(F)}^2 \bigg)^{\frac12}. \label{eq1:lem:bnd_n_sharp}
\end{equation}
\end{Lem}

\bproof 
Let $v\in V\loS+P\upb_k(\calT_h)$ and $w_h\in P\upb_k(\calT_h)$. 
Owing to the definition~\eqref{eq:def_n_sharp} of $n_\sharp$
and the estimate~\eqref{eq:bnd_normal_dual_h} 
from Lemma~\ref{lem:bnd_normal_contrast}, we infer that
\begin{align*}
&|n_\sharp(v,w_h)| \le c\!\sum_{F\in\calFh}\sum_{K\in\calT_F} \!\theta_{K,F}h_K^{d(\frac{1}{2}-\frac1p)}\|\bsigma(v)_{|K}\|_{\bS\upd(K)}  h_F^{-\frac12}\|\jump{w_h}\|_{L^2(F)} \\
&\le c\bigg(\sum_{F\in\calFh}\sum_{K\in\calT_F} \lambda_K^{-\frac12}h_K^{d(\frac{1}{2}-\frac1p)}
\|\bsigma(v)_{|K}\|_{\bL^p(K)}  \mes{\calT_F}^{-\frac12}
\lambda_F^{\frac12}  h_F^{-\frac12}\|\jump{w_h}\|_{L^2(F)} \\
&+ \sum_{F\in\calFh}\sum_{K\in\calT_F} \lambda_K^{-\frac12}h_K^{d(\frac{2+d}{2d}-\frac{1}{q})}
\|\DIV\bsigma(v)_{|K}\|_{L^q(K)} \mes{\calT_F}^{-\frac12}
\lambda_F^{\frac12}  h_F^{-\frac12}\|\jump{w_h}\|_{L^2(F)}\bigg),
\end{align*}
where we used that $\theta_{K,F}\le \theta_{K,F}^{\frac12}$ 
(since $\theta_{K,F}\le 1$), 
$\mes{\calT_F}\lambda_{K}\theta_{K,F}=\lambda_{F}$, the definition 
of $\|\SCAL\|_{\bS\upd(K)}$, and $1+d(\frac12-\frac1q)
=d(\frac{2+d}{2d}-\frac{1}{q})$. 
Owing to the Cauchy--Schwarz inequality, we infer that
$\sum_{F\in\calFh}\sum_{K\in\calT_F} a_K\mes{\calT_F}^{-\frac12}b_F
\le (\sum_{K\in\calT_h}\mes{\calFK}a_K^2)^{\frac12} (\sum_{F\in\calFh} b_F^2)^{\frac12}$,
for all real numbers $\{a_K\}_{K\in\calT_h}$, $\{b_F\}_{F\in\calFh}$, where
we used $\sum_{F\in\calFh}\sum_{K\in\calT_F} = \sum_{K\in\calT_h}\sum_{F\in\calFK}$ for
the term involving the $a_K$'s.
Since $\mes{\calFK}$ is uniformly bounded
($\mes{\calFK}=d+1$ for simplicial meshes), applying this bound to
the two terms composing the above estimate on $|n_\sharp(v,w_h)|$ leads to~\eqref{eq1:lem:bnd_n_sharp}. 
\end{proof}

\bRem[Literature]
Diffusion-dependent averages have been introduced in \citet{Dryja:03}
for discontinuous Galerkin methods and have been analyzed, \eg 
in~\citet{BurZu:06,DryGS:07,DiErG:08,ErStZ:09}. 
\end{Rem}

\section{Applications} \label{Sec:applications}

The goal of this section is to perform a unified error analysis for the
approximation of the model problem~\eqref{model_elliptic_contrast}
with various nonconforming methods: Crouzeix--Raviart finite elements, 
Nitsche's boundary penalty, interior penalty discontinuous Galerkin,
and hybrid high-order methods.
We assume that the exact solution is in
the functional space $V\loS$ defined in~\eqref{eq:def_V_loS_contrast}
with real numbers $p,q$ satisfying~\eqref{eq:p_q_face_to_cell}. 
Our unified analysis hinges on the
dimensionally-consistent seminorm 
\begin{equation}
|v|_{\lambda,p,q}^2 := \|\lambda^{\frac12}\GRAD_hv\|_{\Ldeuxd}^2 + 
|v|_{n_\sharp}^2, \qquad \forall v\in V\loS+P\upb_k(\calT_h),
\end{equation}%
with $|\SCAL|_{n_\sharp}$ defined in~\eqref{eq:seminorm_sharp_n}.
Since $\lambda$ is piecewise constant, we have
\begin{align} 
|v|_{\lambda,p,q}^2 := {}&\sum_{K\in\calT_h} \lambda_K \Big(\|\GRAD v_{|K}\|_{\bL^2(K)}^2
+ h_K^{2d(\frac12-\frac1p)}\|\nabla v_{|K}\|_{\bL^p(K)}^2 \nonumber \\ &+
h_K^{2d(\frac{d+2}{2d}-\frac{1}{q})}\|\Delta v_{|K}\|_{L^q(K)}^2\Big).\label{eq:snorme_lpq}
\end{align}
Invoking inverse inequalities shows that there is $c$, 
uniform \wrt $\hinH$, but depending on $p$ and $q$, \sth
\begin{equation} \label{eq:inv_ineq_pq}
|v_h |_{\lambda,p,q} \le c\,  \|\lambda^{\frac12} \GRAD v_{h}\|_{\bL^2(\Dom)},
\qquad \forall v_h \in P\upb_k(\calT_h).
\end{equation}

\subsection{Abstract approximation result} \label{Sec:abstract_error}
We start by recalling a general approximation result established in
\citep[Lem.~4.4]{Ern_Guermond_cmam_2017}. 
Let $V$ and $W$ be two real Banach spaces. Let $a(\SCAL,\SCAL)$ be a bounded
bilinear form on $V\CROSS W$, and let
$\form(\SCAL)$ be a bounded linear form on $W$, \ie $\form\in W'$.  We
consider the following abstract model problem:
\begin{equation}
\left\{
\begin{array}{l}
\text{Find $u \in V$ such that}\\[2pt]
a(u,w)=\form(w), \quad \forall w \in W,
\end{array}
\right.
\label{eq:gal_weak}
\end{equation}
which we assume to be well-posed in the sense of Hadamard; that is to
say, there is a unique solution and this solution depends continuously
on the data. 

We now formulate a discrete version of the problem~\eqref{eq:gal_weak}
by using the Galerkin method.  We replace the infinite-dimensional
spaces $V$ and $W$ by finite-dimensional spaces $V_h$ and $W_h$ that
are members of sequences of spaces $(V_h)_{\hinH}$, $(W_h)_{\hinH}$
endowed with some approximation properties as $h\to 0$. The norms in
$V_h$ and $W_h$ are denoted by $\|\SCAL\|_{V_h}$ and
$\|\SCAL\|_{W_h}$, respectively. The discrete version of~\eqref{eq:gal_weak} is formulated as
follows:
\begin{equation}
\left\{
\begin{array}{l}
\text{Find $u_h \in V_h$ such that}\\[2pt]
a_h(u_h,w_h)=\form_h(w_h), \quad \forall w_h \in W_h,
\end{array}
\right.
\label{eq:discret_contrast}
\end{equation}
where $a_h(\SCAL,\SCAL)$ is a bounded bilinear form on $V_h\CROSS W_h$ and
$\form_h(\SCAL)$ is a bounded linear form on $W_h$; note that $a_h(\SCAL,\SCAL)$ and
$\form_h(\SCAL)$ possibly differ from $a(\SCAL,\SCAL)$ and $\form(\SCAL)$, respectively.
We henceforth assume that
$\dim(V_h)=\dim(W_h)$ and that
\begin{equation} 
\inf_{0\ne v_h\in V_h} \sup_{0\ne w_h\in W_h}
  \frac{|a_h(v_h,w_h)|}{\| v_h \|_{V_h} \|w_h\|_{W_h}} =:
  \alpha_h>0, \qquad \forall h>0, \label{eq:bnb1hhhh}
\end{equation}
so that the discrete
problem~\eqref{eq:discret_contrast} is well-posed.

We formalize the fact that 
the error analysis 
requires the solution to \eqref{eq:gal_weak} to be slightly more 
regular than just being a member of $V$ by introducing
a functional space $V\loS$ such that $u\in V\loS \subsetneq V$. 
Our setting for the error analysis is therefore as follows:
\begin{equation}\label{eq:Vshs=V+Vh}
u\in V\loS \subsetneq V, \qquad u-u_h \in \Vshs := V\loS+V_h,
\end{equation}
with the norm in $\Vshs$ denoted by $\|{\cdot}\|_{\Vshs}$.
Since $V_h$ is finite-dimensional, we have 
\begin{equation} \label{eq:equiv_Vh_Vshs}
c_{\sharp h}:=\sup_{0\ne v_h\in V_h} \frac{\|v_h\|_{\Vshs} }{\|v_h\|_{V_h}} < \infty.
\end{equation}

We now introduce the consistency error mapping 
$\delta_h : V_h \to W_h':=\calL(W_h;\Real)$
defined for all $v_h\in V_h$  and all $w_h\in W_h$ by setting 
\begin{equation} \label{eq:def_consis_error1} 
\langle \delta_h(v_h),w_h\rangle_{W_h',W_h} := \ell_h(w_h)-a_h(v_h,w_h)
=a_h(u_h-v_h,w_h).
\end{equation}
We further assume that 
\begin{equation}
\omega_{\sharp h} := \sup_{u\in V\loS}\sup_{v_h\in V_h{\setminus}\{u\}} 
\frac{\|\delta_h(v_h)\|_{W_h'}}{\|u-v_h\|_{\Vshs}}< \infty. \label{eq:consistency_abst}
\end{equation}

\begin{Exp}[Conforming setting] \label{exp:simple_setting_error}
Assume conformity, $a_h=a$, and 
$\ell_h=\ell$. Take $V\loS:=V$, so that $\Vshs=V$, and 
take $\|\SCAL\|_{\Vshs}:=\|\SCAL\|_V$. 
The consistency error~\eqref{eq:def_consis_error1} is such that
\[
\langle \delta_h(v_h),w_h\rangle_{W_h',W_h} = \ell(w_h)-a(v_h,w_h) = a(u-v_h,w_h),
\]
where we used that $\ell(w_h)=a(u,w_h)$
(\ie the Galerkin orthogonality property). 
Since $a$ is bounded on $V\CROSS W$, 
\eqref{eq:consistency_abst} holds true with 
$\omega_{\sharp h}=\|a\|$; moreover, $c_{\sharp h}=1$.
\end{Exp}

The main result we are going to invoke later is the following.
\begin{Lem}[Quasi-optimal error estimate] \label{lem:abst_err_est_QO}
If $u\in V\loS$, then
\begin{equation} \label{eq:abst_err_est_QO}
\|u-u_h\|_{\Vshs} \le 
\bigg( 1+c_{\sharp h}\frac{\omega_{\sharp h}}{\alpha_h}\bigg) \inf_{v_h\in V_h} \|u-v_h\|_{\Vshs}.
\end{equation}
\end{Lem}

\bproof
The proof is classical; we sketch it for completeness. 
For all $v_h\in V_h$, we have
\begin{align*}
\|u_h-v_h\|_{\Vshs} & \le c_{\sharp h}\, \|u_h-v_h\|_{V_h} 
\le \frac{c_{\sharp h}}{\alpha_h} \, \sup_{0\ne w_h\in W_h} 
\frac{|a_h(u_h-v_h,w_h)|}{\|w_h\|_{W_h}} \\
&=\frac{c_{\sharp h}}{\alpha_h} \, \|\delta_h(v_h)\|_{W_h'} 
\le \frac{c_{\sharp h}\omega_{\sharp h}}{\alpha_h} \, \|u-v_h\|_{\Vshs}.
\end{align*}
We conclude by using the triangle inequality and
taking the infimum over $v_h\in V_h$.
\end{proof}

When the constants $c_{\sharp h}$ and $\omega_{\sharp h}$ can be bounded 
from above uniformly \wrt $\hinH$, we denote by $c_\sharp$ and $\omega_\sharp$
any constant such that $c_\sharp\ge \sup_{\hinH} c_{\sharp h}$ and 
$\omega_\sharp\ge \sup_{\hinH} \omega_{\sharp h}$.

\subsection{Crouzeix--Raviart approximation}
\label{sec:CR_contrast}

We consider in this section the approximation of the model
problem~\eqref{eq:weak_PDE_contrast} with a homogeneous Dirichlet
condition (for simplicity) using the
Crouzeix--Raviart finite element space
\begin{equation}
P_{1,0}\uCR(\calT_h) := \bset v_h \in P_1\upb(\calT_h)\tq  
\int_F \jump{v_h}_F\dif s = 0,\,\forall F \in \calFh \eset.
\end{equation}
The discrete problem~\eqref{eq:discret_contrast} is formulated 
with $V_h:=P_{1,0}\uCR(\calT_h)$ and the following forms:
\begin{equation} \label{eq:forms_CR}
a_h(v_h,w_h):= \int_\Dom \lambda\GRAD_hv_h\SCAL\GRAD_hw_h\dif x,
\qquad
\ell_h(w_h)=\int_\Dom fw_h\dif x.
\end{equation}
We equip $V_h$ with the norm
$\|v_h\|_{V_h} :=\|\lambda^{\frac12}\GRAD_hv_h\|_{\Ldeuxd}$.
The following result is standard.

\begin{Lem}[Coercivity, well-posedness] \label{lem:stab_CR_contrast} 
The bilinear form
$a_h$ is coercive on $V_h$ with coercivity constant $\alpha=1$, and the discrete
problem~\eqref{eq:discret_contrast} is well-posed.
\end{Lem}

Let $\Vshs:= V\loS + V_h$ be equipped with the norm
$\|v\|_{\Vshs} := |v|_{\lambda,p,q}$ with $|v|_{\lambda,p,q}$ defined
in~\eqref{eq:snorme_lpq} (this is indeed a norm on $\Vshs$ since
$|v|_{\lambda,p,q}=0$ implies that $v$ is piecewise constant and hence
vanishes identically owing to the definition of $V_h$).  
Owing to~\eqref{eq:inv_ineq_pq}, 
there is $c_{\sharp}$, uniform \wrt $\hinH$, but depending on $p$ and
$q$, \sth
$\|v_h\|_{\Vshs} \le c_{\sharp} \|v_h\|_{V_h}$, for all $v_h\in V_h$.  

\begin{Lem}[Consistency/boundedness] \label{lem:consist_CR_contrast}
There is $\omega_\sharp$, uniform \wrt $\hinH$, $\lambda$,
and $u\in V\loS$, but depending on $p$ and $q$, \sth
$\|\delta_h(v_h)\|_{V_h'} \le \omega_\sharp \|u-v_h\|_{\Vshs}$, for
all $v_h\in V_h$.
\end{Lem}

\bproof 
Let $v_h,w_h\in V_h$. Since $V_h\subset P\upb_k(\calT_h)$,
the identity~\eqref{eq:n_sharp_1} implies that
\[
n_\sharp(v_h,w_h) 
= \sum_{F\in\calFh} \int_F \avg{\bsigma(v_h)}_{\theta}\SCAL\bn_F \jump{w_h} \dif s = 0,
\]
because $\avg{\bsigma(v_h)}_{\theta}\SCAL\bn_F$ is constant over $F$.
Moreover, invoking the identity~\eqref{eq:n_sharp_2} 
with $v=u$ and since $f=\DIV\bsigma(u)$, we have
\[
\ell_h(w_h) = n_\sharp(u,w_h)- \int_\Dom \bsigma(u)\SCAL\GRAD_hw_h\dif x.
\]
Combining the two above identities and letting $\eta:=u-v_h$,
we obtain
\begin{align*}
\langle \delta_h(v_h),w_h\rangle_{V_h',V_h} 
&= n_\sharp(u,w_h)  + \int_\Dom \lambda\GRAD_h\eta\SCAL\GRAD_h w_h\dif x
= n_\sharp(\eta,w_h) + \int_\Dom \lambda\GRAD_h\eta\SCAL\GRAD_h w_h\dif x.
\end{align*}%
The first term on the right-hand side is estimated by invoking the
boundedness of $n_\sharp$ (Lemma~\ref{lem:bnd_n_sharp}), the
inequality $\lambda_F\le \min_{K\in\calT_F}\lambda_K$ (see
\eqref{def_of_lambdaF}), and the bound
$\sum_{F\in\calFh} \lambda_F h_F^{-1}\|\jump{w_h}\|_{L^2(F)}^2 \le
c\|w_h\|_{V_h}^2$, which is standard for Crouzeix--Raviart
elements. The second term is estimated by using the Cauchy--Schwarz
inequality.
\end{proof}

\begin{Th}[Error estimate]\label{th:conv_CR_contrast}
Let $u$ solve~\eqref{eq:weak_PDE_contrast} and $u_h$ 
solve~\eqref{eq:discret_contrast} with $a_h$ and $\ell_h$ defined
in~\eqref{eq:forms_CR}.
Assume that there is $r>0$ \sth $u\in H^{1+r}(\Dom)$. 
There is $c$, uniform \wrt $\hinH$, $\lambda$,
and $u\in H^{1+r}(\Dom)$, 
but depending on $r$, \sth
the following quasi-optimal error estimate holds true:
\begin{equation} \label{eq:QO_CR_contrast}
\|u-u_h\|_{\Vshs} \le c\, \inf_{v_h\in V_h} \|u-v_h\|_{\Vshs}.
\end{equation}
Moreover, letting $t:=min(1,r)$, where $1=k$ is the 
degree of the Crouzeix--Raviart finite element, we have
\begin{equation} \label{eq:rate_CR_contrast} 
\|u-u_h\|_{\Vshs} \le c\,
\bigg(\sum_{K\in\calT_h}\!\lambda_Kh_K^{2t}|u|_{H^{1+t}(K)}^2
+\lambda_K^{-1} h_K^{2d(\frac{d+2}{2d}-\frac{1}{q})}\|f\|_{L^q(K)}^2\bigg)^{\frac{1}{2}}\!.
\end{equation} 
\end{Th}
\bproof The error estimate~\eqref{eq:QO_CR_contrast} follows from
Lemma~\ref{lem:abst_err_est_QO} combined with stability
(Lemma~\ref{lem:stab_CR_contrast}) and consistency/boundedness
(Lemma~\ref{lem:consist_CR_contrast}).  We now bound the infimum
in~\eqref{eq:QO_CR_contrast} by considering $\eta:=u-\inter_h\uCR(u)$,
where $\inter_h\uCR$ is the Crouzeix--Raviart interpolation
operator using averages over the faces as degrees of freedom. 
It is a standard approximation result 
that there is $c$, uniform \wrt $u\in H^{1+t}(K)$, $t\ge 0$, and $\hinH$, 
\sth $\|\GRAD\eta_{|K}\|_{\bL^2(K)} \le ch_K^t|u|_{H^{1+t}(K)}$ for all $K\in\calT_h$.
Moreover, invoking the embedding
$\bH^t(\wK)\hookrightarrow \bL^p(\wK)$ and classical results on the 
transformation of Sobolev norms by the geometric mapping, we obtain the bound
\begin{equation} \label{eq:inv_ineq_Lp_Hr_contrast}
h_K^{d(\frac12-\frac{1}{p})}\|\GRAD\eta_{|K}\|_{\bL^p(K)} \le c\,
\big(\|\GRAD\eta_{|K}\|_{\bL^2(K)} + h_{K}^t
|\GRAD\eta_{|K}|_{\bH^t(K)}\big).
\end{equation}
Observing that $|\GRAD \eta_{|K}|_{H^{t}(K)}=|u|_{H^{1+t}(K)}$ since
$\inter_h\uCR(u)$ is affine on $K$ and using again the approximation properties
of $\inter_h\uCR$, we infer that
$h_K^{d(\frac12-\frac1p)}\|\GRAD\eta_{|K}\|_{\bL^p(K)} \le c\,
h_K^t|u|_{H^{1+t}(K)}$.
Finally, we have $\Delta \eta_{|K} = \lambda_K^{-1}f$ in $K$.  
\end{proof}

\bRem[Convergence]
The rightmost term in \eqref{eq:rate_CR_contrast} 
converges as $O(h)$ when $q=2$. Moreover, convergence is lost when
$q\le\frac{2d}{d+2}$, which is somewhat natural since in this case the
linear form $w\mapsto \int_\Dom fw \dif x$ is no longer bounded on
$\Hun$.
\end{Rem}

\bRem[Weights] Although the weights introduced
in~\eqref{eq:def_weight_theta} are not explicitly used in the Crouzeix--Raviart
discretization, they play a role in the error analysis. More
precisely, we used the boundedness of 
the bilinear form $n_\sharp$ together with 
$\lambda_F\le \min_{K\in\calT_F}\lambda_K$ in the proof of
Lemma~\ref{lem:consist_CR_contrast}. 
The present approach is somewhat more
general than that in \cite{Li_Shipeng_2013} since it delivers error estimates that are robust
with respect to the diffusivity contrast. The trimming operator 
invoked in \citep[Eq. (5)--(7)]{Li_Shipeng_2013} cannot account for the 
diffusivity contrast.
\end{Rem}

\subsection{Nitsche's boundary penalty method}
\label{sec:Nitsche_contrast}

We consider in this section the approximation of the model
problem~\eqref{model_elliptic_contrast} by means of Nitsche's boundary
penalty method.  Now we set
\begin{equation}
V_h:=P\upg_k(\calT_h) := \bset v_h\in
P\upb_k(\calT_h) \tq \jump{v_h}_F = 0,\;\forall F\in\calFhi\eset, \qquad
k\ge1,
\end{equation}
\ie
$V_h$ is $H^1$-conforming
The discrete problem~\eqref{eq:discret_contrast} is formulated with
$V_h:=P\upg_k(\calT_h)$ and the following forms:
\begin{subequations} \label{eq:forms_Nitsche} \begin{align} 
a_h(v_h,w_h) &:= a(v_h,w_h) + \sum_{F\in \calFhb} 
\int_F \bigg(\bsigma(v_h)\SCAL \bn+\varpi_0\frac{\lambda_{K_l}}{h_F}v_h\bigg)w_h\dif s,\\
\form_h(w_h) &:= \ell(w_h) + \sum_{F\in \calFhb} 
\varpi_0\frac{\lambda_{K_l}}{h_F} \int_F  gw_h\dif s,
\end{align}\end{subequations}
where the exact forms $a$ and $\ell$ are defined in~\eqref{eq:a_l_contrast},
$K_l$ is the
unique mesh cell \sth $F=\partial K_l\cap \front$, and the user-specified
penalty parameter $\varpi_0$ is yet to be chosen large enough. It is
possible to add a symmetrizing term to the discrete bilinear form $a_h$. 

We equip $V_h$ with the norm
$\|v_h\|_{V_h}^2 := \|\lambda^{\frac12}\GRAD v_h\|_{\Ldeuxd}^2 + \bndsn{v_h}^2$ with
$\bndsn{v_h}^2 := \sum_{F\in\calFhb} \frac{\lambda_{K_l}}{h_F}\|v_h\|_{L^2(F)}^2$.
Owing to the shape-regularity of the mesh sequence, there is
$c_I$, uniform \wrt $\hinH$ \sth
\begin{equation} \label{inv_trace_ineq_Nitsche}
\|v_h\|_{L^2(F)}\le c_Ih_F^{-\frac12} \|v_h\|_{L^2(K_l)},
\end{equation}
for all $v_h\in V_h$ and all $F\in \calFhb$. 
Let $n_\partial$ denote the 
maximum number of boundary faces that a mesh cell 
can have ($n_\partial\le d$ for simplicial meshes).
The proof of the following result uses
standard arguments.

\begin{Lem}[Coercivity, well-posedness] \label{lem:stab_Nitsche_contrast} 
Assume that
the penalty parameter satisfies $\varpi_0>\frac14n_\partial c_I^2$.
Then, $a_h$ is coercive on $V_h$ with constant
$\alpha:=\frac{\varpi_0-\frac14n_\partial c_I^2}{1+\varpi_0}>0$, and the
discrete problem~\eqref{eq:discret_contrast} is well-posed.
\end{Lem}


Let $\Vshs:=V\loS+V_h$. We equip the space $\Vshs$ 
with the norm $\|v\|_{\Vshs}^2
:= |v|_{\lambda,p,q}^2 + \bndsn{v}^2$ with 
\begin{align} 
|v|_{\lambda,p,q}^2 := {}&\sum_{K\in\calT_h} \lambda_K \|\GRAD v_{|K}\|_{\bL^2(K)}^2 \nonumber \\&
+ \sum_{K\in\calThbb} \lambda_K\Big(
h_K^{2d(\frac12-\frac1p)}\|\nabla v_{|K}\|_{\bL^p(K)}^2 +
h_K^{2d(\frac{d+2}{2d}-\frac{1}{q})}\|\Delta v_{|K}\|_{L^q(K)}^2\Big),\label{eq:snorme_lpq_Nitsche}
\end{align}
where $\calThbb$ is the collection of the mesh cells having
at least one boundary face, and $\bndsn{v}^2 = \sum_{F\in\calFhb} 
\frac{\lambda_{K_l}}{h_F}\|v\|_{L^2(F)}^2$.
Owing to~\eqref{eq:inv_ineq_pq}, there is $c_{\sharp}$, uniform \wrt $\hinH$, 
but depending on $p$ and $q$, \sth $\|v_h\|_{\Vshs} \le c_{\sharp} 
\|v_h\|_{V_h}$, for all $v_h\in V_h$.

\begin{Lem}[Consistency/boundedness] \label{lem:consist_Nitsche_contrast}
There is $\omega_\sharp$, uniform \wrt $\hinH$, $\lambda$,
and $u\in V\loS$, but depending on $p$ and $q$, \sth
$\|\delta_h(v_h)\|_{V_h'} \le \omega_\sharp \|u-v_h\|_{\Vshs}$, for
all $v_h \in V_h$.
\end{Lem}

\bproof Let $v_h,w_h\in V_h$.  Using the identity~\eqref{eq:n_sharp_1}
for $n_\sharp$, $\jump{w_h}_F=0$ for all $F\in\calFhi$ (since $V_h$
is $H^1$-conforming), and the definition of the weights at the boundary faces, we infer that
$n_\sharp(v_h,w_h)=\sum_{F\in\calFhb}\int_F \bsigma(v_h)\SCAL\bn w_h
\dif s$. Hence,
$a_h(v_h,w_h) = a(v_h,w_h) + n_\sharp(v_h,w_h) 
+ \sum_{F\in\calFhb} \varpi_0\frac{\lambda_{K_l}}{h_F} \int_F v_hw_h\dif s$.
Therefore, invoking the identity~\eqref{eq:n_sharp_2} for the exact
solution $u$ and observing that $f=\DIV\bsigma(u)$,
we infer the important identity $\int_\Dom f w_h\dif x =  
a (u,w_h) + n_\sharp (u,w_h)$.
Then, recalling that $\gamma\upg(u)=g$,
and letting $\eta:=u-v_h$, we obtain
\[
\langle \delta_h(v_h),w_h\rangle_{V_h',V_h} 
= n_\sharp(\eta,w_h) + a(\eta,w_h) + \sum_{F\in\calFhb} \varpi_0\frac{\lambda_{K_l}}{h_F} \int_F \eta w_h\dif s.
\]
We conclude by using the boundedness
of $n_\sharp$ from Lemma~\ref{lem:bnd_n_sharp} and the Cauchy--Schwarz inequality.
\end{proof}

\begin{Th}[Error estimate] \label{th:conv_Nitsche_contrasted} 
Let $u$ solve~\eqref{model_elliptic_contrast} and $u_h$
solve~\eqref{eq:discret_contrast} with $a_h$ and $\ell_h$ 
defined in~\eqref{eq:forms_Nitsche} and penalty parameter
$\varpi_0>\frac14n_\partial c_I^2$. Assume that there is $r>0$ \sth
$u\in H^{1+r}(\Dom)$.  There is $c$, uniform with respect to
$\hinH$, $\lambda$, and $u\in H^{1+r}(\Dom)$, but depending on
$r$, \sth the following quasi-optimal error estimate holds true:
\begin{equation} \label{eq:QO_Nitsche_contrast}
\|u-u_h\|_{\Vshs} \le c\, \inf_{v_h\in V_h} \|u-v_h\|_{\Vshs}.
\end{equation}
Moreover, letting $t:=min(r,k)$, $\chi_t=1$ if $t\le 1$ and $\chi_t=0$ if $t>1$, we have
\begin{equation} \label{eq:rate_Nitsche_contrast}
\|u-u_h\|_{\Vshs} \le c\, \bigg(\sum_{K\in\calT_h} \!\!\lambda_Kh_K^{2t}|u|_{H^{1+t}(\check\calT_K)}^2
+ \frac{\chi_t}{\lambda_K} h_K^{2d(\frac{d+2}{2d}-\frac{1}{q})}\|f\|_{L^q(K)}^2
\bigg)^{\frac12}\!,
\end{equation}
where $\check\calT_K$ is the collection of the mesh cells having at
least a common vertex with $K$. The broken Sobolev norm 
$|\SCAL|_{H^{1+t}(\check\calT_K)}$ can
be replaced by $|\SCAL|_{H^{1+t}(K)}$ if $1+t>\frac{d}{2}$.
\end{Th}

\bproof 
The error estimate~\eqref{eq:QO_Nitsche_contrast} follows from 
Lemma~\ref{lem:abst_err_est_QO} combined with 
stability (Lemma~\ref{lem:stab_Nitsche_contrast}) and
consistency/boundedness (Lemma~\ref{lem:consist_Nitsche_contrast}).
We now bound the 
infimum in~\eqref{eq:QO_Nitsche_contrast} by using
$\eta:=u-\inter_h\upgav(u)$, where $\inter_h\upgav$ is the 
quasi-interpolation operator introduced 
in~\citep[\S5]{Ern_Guermond_M2AN_2017}. We take 
the polynomial degree of $\inter_h\upgav$ to be
$\ell:=\lceil t\rceil$, where $\lceil t\rceil$ denotes
the smallest integer $n\in\polN$ \sth $n\ge t$. Notice that $\ell\ge 1$ because $r>0$ and $k\ge 1$, and
$\ell\le k$ because 
$t\le k$; hence, $\inter_h\upgav(u)\in V_h$.
We need to bound all the terms composing the norm $\|\eta\|_{\Vshs}$. Owing to
\citep[Thm.~5.2]{Ern_Guermond_M2AN_2017} with $m=1$, we have
$\|\GRAD \eta\|_{\bL^2(K)} \le ch_K^{t}|u|_{H^{1+t}(\check\calT_K)}$
for all $K\in\calT_h$. Moreover, we have
$h_F^{-\frac12}\|\eta\|_{L^2(F)} \le
ch_{K_l}^{t}|u|_{H^{1+t}(\check\calT_{K_l})}$ for all $F\in\calFhb$.
It remains to estimate
$h_K^{d(\frac12-\frac1p)}\|\GRAD\eta_{|K}\|_{\bL^p(K)}$ and
$h_K^{d(\frac{d+2}{2d}-\frac{1}{q})}\|\Delta\eta_{|K}\|_{L^q(K)}$ for all
$K\in\calThbb$.  Using~\eqref{eq:inv_ineq_Lp_Hr_contrast}, the above
bound on $\|\GRAD \eta\|_{\bL^2(K)}$, and
$|\GRAD\eta|_{\bH^{t}(K)} =|\GRAD u|_{\bH^{t}(K)}=|u|_{H^{1+t}(K)}$
since $\ell<1+t$, we infer that
$h_K^{d(\frac12-\frac1p)}\|\GRAD\eta\|_{\bL^p(K)} \le c\,
h_K^{t}|u|_{H^{1+t}(\check\calT_K)}$.
Moreover, if $t\le1$, we have $\ell=1$ so that
$\|\Delta \eta_{|K}\|_{L^{q}(K)}=\|\Delta u\|_{L^{q}(K)} =
\lambda_K^{-1}\|f\|_{L^{q}(K)}$.
Instead, if $t>1$, we infer that $r>1$ so that we can set $q=2$
(recall that $f_{|\Dom_i} = \lambda_{|\Dom_i}(\Delta u)_{\Dom_i}$ for
all $i\in\intset{1}{M}$, and $u\in H^2(\Dom)$ if $r\ge1$), and we
estimate $\|\Delta\eta_{|K}\|_{L^2(K)}$ using
\citep[Thm.~5.2]{Ern_Guermond_M2AN_2017} with $m=2$.  Finally, if
$1+t>\frac{d}{2}$, we can use the canonical Lagrange interpolation operator
$\inter_h\upg$ instead of $\inter_h\upgav$, and this allows us to
replace $|\SCAL|_{H^{1+t}(\check\calT_K)}$ by $|\SCAL|_{H^{1+t}(K)}$ 
in~\eqref{eq:rate_Nitsche_contrast}. 
\end{proof}


\subsection{Discontinuous Galerkin}
\label{sec:dG_contrast}

We consider in this section the
approximation of the model problem~\eqref{model_elliptic_contrast}
by means of the symmetric interior penalty discontinuous Galerkin
method.
The discrete problem~\eqref{eq:discret_contrast} is formulated with 
$V_h:=P_k\upb(\calT_h)$, $k\ge1$, the bilinear forms
\begin{subequations} \label{eq:forms_IPDG} \begin{align} 
a_h(v_h,w_h) :={}& \int_\Dom \lambda \GRAD_h v_h\SCAL \GRAD_h w_h\dif x 
+ \sum_{F\in\calFh} \int_F \avg{\bsigma(v_h)}_\theta\SCAL \bn_F\jump{w_h}\dif s \nonumber \\
&+\sum_{F\in\calFh} \int_F \jump{v_h}\avg{\bsigma(w_h)}_\theta\SCAL \bn_F\dif s 
+ \sum_{F\in\calFh} \varpi_0\frac{\lambda_{F}}{h_F}\int_F \jump{v_h}\jump{w_h}\dif s, \\
\form_h(w_h) := {}&\ell(w_h) + \sum_{F\in \calFhb}
\varpi_0\frac{\lambda_{K_l}}{h_F} \int_F g w_h\dif s,
\end{align}\end{subequations}
where $\ell$ is
defined in~\eqref{eq:a_l_contrast}, $\lambda_F$ in~\eqref{def_of_lambdaF},
and the user-specified penalty
parameter $\varpi_0$ is yet to be chosen large enough.  We equip $V_h$
with the norm
$\|v_h\|_{V_h}^2 := \|\lambda^{\frac12}\GRAD_h v_h\|_{\Ldeuxd}^2 +
\jumpsn{v_h}^2$ with
$\jumpsn{v_h}^2 := \sum_{F\in\calFh}
\frac{\lambda_{F}}{h_F}\|\jump{v_h}\|_{L^2(F)}^2$.  Recall the
discrete trace inequality~\eqref{inv_trace_ineq_Nitsche} and 
recall that $n_\partial$ denotes the maximum number of faces that a mesh cell can
have ($n_\partial\le d+1$ for simplicial meshes). The proof of the following 
result uses standard arguments.

\begin{Lem}[Coercivity, well-posedness] \label{lem:coer_dg_contrast} 
Assume that the penalty parameter satisfies
$\varpi_0>n_\partial c_I^2$. Then, $a_h$ is coercive on $V_h$ with constant
$\alpha:=\frac{\varpi_0-n_\partial c_I^2}{1+\varpi_0}>0$, and the discrete 
problem~\eqref{eq:discret_contrast} is well-posed. 
\end{Lem}


Let $\Vshs:= V\loS + V_h$. We equip the space $\Vshs$ 
with the norm 
$\|v\|_{\Vshs}^2 := |v|_{\lambda,p,q}^2+\jumpsn{v}^2$ 
with $|v|_{\lambda,p,q}$ defined in~\eqref{eq:snorme_lpq} and
$\jumpsn{v}^2:= \sum_{F\in\calFh} \frac{\lambda_{F}}{h_F}
\|\jump{v}\|_{L^2(F)}^2$.
Owing to~\eqref{eq:inv_ineq_pq}, there is $c_{\sharp}$, 
uniform \wrt $\hinH$, but depending on $p$ and $q$, 
\sth $\|v_h\|_{\Vshs} \le c_{\sharp} \|v_h\|_{V_h}$, for all $v_h\in V_h$.

\begin{Lem}[Consistency/boundedness] \label{lem:consist_dG_contrast}
There is $\omega_\sharp$, uniform \wrt $\hinH$, $\lambda$, and 
$u\in V\loS$, but depending on $p$ and $q$, \sth
$\|\delta_h(v_h)\|_{V_h'} \le \omega_\sharp \|u-v_h\|_{\Vshs}$, 
for all $v_h \in V_h$.
\end{Lem}

\bproof Let $v_h,w_h\in V_h$.  Owing to~\eqref{eq:n_sharp_2} and since $f=\DIV\bsigma(u)$, 
we infer that $\int_\Dom f w_h \dif x = \sum_{K\in\calT_h}a_K(u,w_h) 
+ n_\sharp(u,w_h)$ with $a_K(u,w_h):=-(\bsigma(u),\GRAD_h w_h)_{\bL^2(K)}$. 
Using the identity~\eqref{eq:n_sharp_1}, we obtain
\begin{align*}
\ell_h(w_h)={}&n_\sharp(u,w_h) - \int_\Dom \bsigma(u)\SCAL \GRAD_h w_h\dif x 
+ \sum_{F\in\calFhb} \varpi_0\frac{\lambda_F}{h_F}\int_F g w_h \dif s,\\
a_h(v_h,w_h) ={}& \int_\Dom-\bsigma(v_h)\SCAL \GRAD_h w_h\dif x +n_\sharp(v_h,w_h) \\
&-\sum_{F\in\calFh} \int_F \jump{v_h}\avg{\bsigma(w_h)}_\theta\SCAL\bn_F\dif s 
+\sum_{F\in\calFh} \varpi_0\frac{\lambda_F}{h_F}\int_F \jump{v_h} \jump{w_h} \dif s.
\end{align*}
Then setting $\eta := u-v_h$ and using that $\jump{u}_F=0$ for all $F\in\calFhi$ and $\jump{u}_F=g$ for all $F\in\calFhb$, we obtain the following representation of the consistency linear form $\delta_h(v_h)$:
\begin{align*}
\langle \delta_h(v_h),w_h\rangle_{V_h',V_h} = {}& n_\sharp(\eta,w_h) + \int_\Dom\lambda\GRAD\eta\SCAL \GRAD_h w_h\dif x  \\
&-\sum_{F\in\calFh} \int_F \jump{\eta}\avg{\bsigma(w_h)}_\theta\SCAL\bn_F\dif s + \sum_{F\in\calFh} \!\varpi_0\frac{\lambda_F}{h_F}\int_F \jump{\eta} \jump{w_h} \dif s.
\end{align*}
Bounding the second, third and fourth terms uses standard arguments
(see, \eg \citep{DiPEr:12}), whereas we invoke
the boundedness estimate on
$n_\sharp$ from Lemma~\ref{lem:bnd_n_sharp} for the first term.
\end{proof}

\begin{Th}[Error estimate] \label{th:conv_dG_contrasted}
Let $u$ solve~\eqref{model_elliptic_contrast} and $u_h$ 
solve~\eqref{eq:discret_contrast} with $a_h$ and $\ell_h$ 
defined in~\eqref{eq:forms_IPDG} and penalty parameter 
$\varpi_0>n_\partial c_I^2$. Assume that there is $r>0$ \sth $u\in H^{1+r}(\Dom)$. 
There is $c$, uniform with respect to $\hinH$, $\lambda$, 
and $u\in H^{1+r}(\Dom)$, 
but depending on $r$, \sth
the following quasi-optimal error estimate holds true:
\begin{equation} \label{eq:QO_dG_contrast}
\|u-u_h\|_{\Vshs} \le c\, \inf_{v_h\in V_h} \|u-v_h\|_{\Vshs}.
\end{equation}
Moreover, letting $t:=min(r,k)$, $\chi_t=1$ if $t\le1$ and $\chi_t=0$ if $t>1$, we have
\begin{equation} \label{eq:rate_dG_contrast}
\|u-u_h\|_{\Vshs} \le c\,\bigg(\sum_{K\in\calT_h} \!\!\lambda_Kh_K^{2t}|u|_{H^{1+t}(K)}^2
+ \frac{\chi_t}{\lambda_K} 
h_K^{2d(\frac{d+2}{2d}-\frac{1}{q})}\|f\|_{L^q(K)}^2\bigg)^{\frac12}\!.
\end{equation} 
\end{Th}
\bproof 
We proceed as in the proof of
Theorem~\ref{th:conv_Nitsche_contrasted}, where we now use the
$L^1$-stable interpolation operator
$\inter_h^\sharp:L^1(\Dom)\to P\upb_k(\calT_h)$
from~\citep[\S3]{Ern_Guermond_M2AN_2017}  to estimate the best approximation
error.  
\end{proof}

\subsection{Hybrid high-order methods}
\label{sec:HHO_contrast}

We consider in this section the approximation of the model
problem~\eqref{model_elliptic_contrast} with a homogeneous Dirichlet
condition (for simplicity) by means of the hybrid
high-order (HHO) method introduced in \citep{DiPEr:15,DiPEL:14}.  We
consider the discrete product space
$\hat V_{h,0}^k:=V_{\calT_h}^k \times V_{\calFh}^k$ with $k\ge0$,
where
\begin{subequations}\begin{align}
V_{\calT_h}^k &:= \{v_{\calT_h}\in L^2(\Dom)\tq v_K:=v_{\calT_h|K} \in V_K^k,\, 
\forall K\in\calT_h\}, \\
V_{\calFh}^k &:= \{ v_{\calFh}\in 
L^2(\calFh)\tq v_{\partial K}:=v_{\calFh|\partial K} 
\in V_{\partial K}^k, \, \forall K\in\calT_h;\; v_{\calFh|\calFhb}=0\},
\end{align}\end{subequations}
with $V_K^k := \polP_{k,d}$ and $V_{\partial K}^k := 
\{ \theta\in L^2(\partial K)\tq 
\theta\circ \bT_{K|\bT_K^{-1}(F)} \in \polP_{k,d-1}, \, \forall F\in\calF_K\}$. 
Thus, for any pair $\hat v_h := (v_{\calT_h},v_{\calFh})\in \hat V_{h,0}^k$,
$v_{\calT_h}$ a collection of cell polynomials
of degree at most $k$, and $v_{\calFh}$
is a collection of face polynomials
of degree at most $k$ which are single-valued at the mesh interfaces and
vanish at the boundary faces (so as to enforce strongly the homogeneous
Dirichlet condition). 
We use the notation
$\hat v_K:=(v_K,v_{\partial K})\in \hat V_K^k:=V_K^k\times V_{\partial K}^k$ 
for all $K\in\calT_h$. 
We equip the local space $\hat V_K^k$ with the $H^1$-like seminorm
\begin{equation} \label{eq:def_norme_locale_HHO}
|\hat v_K|_{\hat V_K^k}^2 := \|\GRAD v_K\|_{\bL^2(K)}^2
+ \|h_{\partial K}^{-\frac12}(v_K-v_{\partial K})\|_{L^2(\partial K)}^2,
\quad \forall \hat v_K=(v_K,v_{\partial K}) \in \hat V_K^k,
\end{equation}
and the global space
$\hat V_{h,0}^k$ with the norm
\begin{equation} \label{eq:norme_HHO_contraste}
  \|\hat v_h\|_{\hat V_{h,0}^k}^2:=\sum_{K\in\calT_h} \lambda_K
  |\hat v_K|_{\hat V_K^k}^2. 
\end{equation}

We introduce locally in each mesh cell $K\in\calT_h$ 
a reconstruction operator and a stabilization operator.
The reconstruction operator
$\opRec_K^{k+1} : \hat V_K^k \to \polP_{k+1,d}$ is defined such that, for any 
pair $\hat v_K=(v_K,v_{\partial K})\in\hat V_K^k$, the polynomial function
$\opRec_K^{k+1}(\hat v_K) \in \polP_{k+1,d}$ solves
\begin{align}
(\GRAD \opRec_K^{k+1}(\hat v_K),\GRAD q)_{\bL^2(K)} &:=
- (v_K,\Delta q)_{L^2(K)} + (v_{\partial K},\bn_K\SCAL\GRAD q)_{L^2(\partial K)}, \label{eq:def_Rec_HHO}
\end{align} 
for all $q\in \polP_{k+1,d}$, with the mean-value condition
$\int_K (\opRec_K^{k+1}(\hat v_K)-v_K)\dif x = 0$.  
This local Neumann problem makes
sense since the right-hand side of~\eqref{eq:def_Rec_HHO} vanishes
when the test function $q$ is constant. 
The stabilization operator
$\opSt_{\partial K}^k : \hat V_K^k \to V_{\partial K}^k$ is defined \sth
for any pair $\hat v_K=(v_K,v_{\partial K}) \in \hat V_K^k$,
\begin{equation} \label{eq:def_SKk_HHO}
\opSt_{\partial K}^k(\hat v_K) := \Pi_{\partial K}^k\left(v_{K|\partial K}
-v_{\partial K}+((I-\Pi_K^k)\opRec_K^{k+1}(\hat v_K))_{|\partial K}\right),
\end{equation}
where $I$
is the identity, $\Pi_{\partial K}^k: L^2(\partial K) \to V_{\partial K}^k$ is the $L^2$-orthogonal
projection onto $V_{\partial K}^k$ and $\Pi_{K}^k: L^2(K) \to V_{K}^k$ is the $L^2$-orthogonal
projection onto $V_K^k$. Elementary algebra shows that the 
stabilization operator can be rewritten as
\begin{equation}\label{eq:def_SKk_HHO_bis}
\opSt_{\partial K}^k(\hat v_K) = \Pi_{\partial K}^k\left(\delta_{\partial K}-((I-\Pi_K^k)\opRec_K^{k+1}(0,\delta_{\partial K}))_{|\partial K}\right),
\end{equation}
with $\delta_{\partial K}:= v_{K|\partial K}-v_{\partial K}$ is a measure of the 
discrepancy between the trace of the cell unknown and the face unknown.

We now introduce the local bilinear form $\hat a_K$ 
on $\hat V_K^k\times \hat V_K^k$ \sth 
\begin{multline}
\hat a_K(\hat v_K,\hat w_K) := 
(\GRAD\opRec_K^{k+1}(\hat v_K),\GRAD \opRec_K^{k+1}(\hat w_K))_{\Ldeuxd} \\
+ (h_{\partial K}^{-1} \opSt_{\partial K}^k(\hat v_K),
\opSt_{\partial K}^k(\hat w_K))_{L^2(\partial K)},
\end{multline} 
where $h_{\partial K}$ is the piecewise constant function on
$\partial K$ \sth $h_{\partial K|F}:=h_F$ for all $F\in\calF_K$.
Then we set 
\begin{equation}
\hat a_h(\hat v_h,\hat w_h) := 
\sum_{K\in\calT_h} \lambda_K 
\hat a_K(\hat v_K,\hat w_K), \qquad
\hat \ell_h(\hat w_h) := \sum_{K\in\calT_h} (f,w_K)_{L^2(K)}.
\label{eq:def_hat_a_HHO}
\end{equation} 
The discrete problem is finally formulated as follows: Find
$\hat u_h \in \hat V_{h,0}^k$ \sth
\begin{equation}
\hat a_h(\hat u_h,\hat w_h) = \hat \ell_h(\hat w_h), 
\quad \forall \hat w_h \in \hat V_{h,0}^k.
\label{eq:weak_HHO_contrast}
\end{equation}

Notice that HHO methods are somewhat simpler than dG methods 
when it comes to solving problems with contrasted coefficients.
For HHO methods, one assembles cellwise the local
bilinear forms $\hat a_K$ weighted by the local diffusion coefficient 
$\lambda_K$, whereas, for dG methods one has to invoke interface-based
values of the diffusion coefficient to construct the penalty term.

The following result is proved in \citep{DiPEr:15,DiPEL:14}.

\begin{Lem}[Stability, boundedness, well-posedness] \label{lem:stab_HHO} 
There are $0<\alpha\le \omega$, uniform \wrt $\hinH$, such that
\begin{equation*} 
\alpha \, |\hat v_K|_{\hat V_K^k}^2 \le \|\GRAD\opRec_K^{k+1}(\hat v_K)\|_{\bL^2(K)}^2 +  
\|h_{\partial K}^{-\frac12}\opSt_{\partial K}^k(\hat v_K)\|_{L^2(\partial K)}^2
= \hat a_K(\hat v_K,\hat v_K)
\le \omega \, |\hat v_K|_{\hat V_K^k}^2,
\end{equation*}
for all $\hat v_K\in \hat V_K$ and all $K\in\calT_h$, and the discrete 
problem~\eqref{eq:weak_HHO_contrast}
is well-posed.
\end{Lem}

The two key tools in the error analysis of HHO methods are a local
reduction operator and the local elliptic projection. For
all $K\in\calT_h$, the local reduction operator
$\hat\inter_K^k:H^1(K)\to \hat V_K^k$ is defined by
$\hat\inter_K^k(v):=(\Pi_K^k(v), \Pi_{\partial
  K}^k(\gamma\upg_{\partial K}(v)))\in\hat V_K^k$, for all
$v\in H^1(K)$.
The local elliptic projection
$\calE_K^{k+1} : H^1(K) \to \polP_{k+1,d}$ is \sth
$(\GRAD(\calE_K^{k+1}(v)-v),\GRAD q)_{\bL^2(K)} = 0$, for all
$q\in \polP_{k+1,d}$, and $(\calE_K^{k+1}(v)-v,1)_{L^2(K)}=0$.  
The following result is established in \citep{DiPEr:15,DiPEL:14}.

\begin{Lem}[Polynomial invariance] \label{lem:ell_proj_HHO} 
The following holds true: 
\begin{subequations} \label{eq:opRst_opSt_HHO} \begin{align}
\label{eq:ell_proj_HHO}
\opRec_K^{k+1}\circ\hat\inter_K^k &= \calE_K^{k+1}, \\
\label{eq:opSt_inter}
\opSt_{\partial K}^k\circ\hat\inter_K^k 
&= (\gamma\upg_{\partial K}\circ\Pi_K^k
-\Pi_{\partial K}^k\circ\gamma\upg_{\partial K}) \circ  (I-\calE_K^{k+1}).
\end{align}\end{subequations}
In particular, 
$\opRec_K^{k+1}(\hat\inter_K^k(p)) = p$ and $\opSt_{\partial K}^k(\hat\inter_K^k(p)) = 0$ for all $p\in \polP_{k+1,d}$.
\end{Lem}


Recalling the duality pairing $\langle\cdot,\cdot\rangle_F$
defined in \eqref{eq:def_comp_normale}, the generalization of the
bilinear form $n_\sharp$ in the context of HHO methods is the
bilinear form defined on
$(V\loS+P\upb_{k+1}(\calT_h))\times \hat V_{h,0}^k$ that acts as
follows:
\begin{equation}
n_\sharp(v,\hat w_h) := \sum_{K\in\calT_h} \sum_{F\in\calFK} \langle (\bsigma(v)\SCAL\bn_K)_{|F},(w_K-w_{\partial K})_{|F}\rangle_F.
\end{equation}

\begin{Lem}[Identities and boundedness for $n_\sharp$] \label{lem:identity_nsharp_hho}
The following holds true for all $\hat w_h\in \hat V_{h,0}^k$, 
all $v_h\in P\upb_{k+1}(\calT_h)$ and all $v\in V\loS$:
\begin{subequations} \label{eq:magic_formula_hho} 
\begin{align}
n_\sharp(v_h,\hat w_h) &= \! \sum_{K\in\calT_h} \! \int_K \lambda_K\GRAD v_{h|K}\SCAL \GRAD(\opRec_K^{k+1}(\hat w_K)-w_K)\dif x, \label{eq:magic_formula1_HHO} \\
n_\sharp(v,\hat w_h) &= \! \sum_{K\in\calT_h} \int_K \bigg( \bsigma(v)\SCAL\GRAD w_K + (\DIV\bsigma(v))w_K\bigg)\dif x. \label{eq:magic_formula2_HHO}
\end{align}\end{subequations}
Moreover, there is $c$, uniform \wrt $\hinH$ and $\lambda$, but
depending on $p$ and $q$, \sth the following holds true
for all $v\in V\loS+P\upb_{k+1}(\calT_h)$ and all
$\hat w_h\in \hat V_{h,0}^k$:
\begin{equation} \label{eq:bnd_nsharp_HHO}
|n_\sharp(v,\hat w_h)| \le c\, |v|_{n_\sharp} 
\bigg( \sum_{K\in\calT_h} \lambda_Kh_K^{-1}\|w_K-w_{\partial K}\|_{L^2(\partial K)}^2 
\bigg)^{\frac12},
\end{equation}
with the $|\SCAL|_{n_\sharp}$-seminorm defined in~\eqref{eq:seminorm_sharp_n}.
\end{Lem}

\bproof \textup{(i)} We first prove \eqref{eq:magic_formula1_HHO}.
Let $v_h\in P\upb_{k+1}(\calT_h)$ and $\hat w_h\in \hat V_{h,0}^k$. 
Since the restriction of $\bsigma(v_h)$ to each mesh cell is smooth
and since the trace on $\partial K$ of the face-to-cell lifting 
operator $L_F^K$ is nonzero only on $F$, for all $F\in\calFK$, we have
\begin{align*}
&\langle (\bsigma(v_h)\SCAL\bn_K)_{|F},(w_K-w_{\partial K})_{|F}\rangle_F \\
&=\int_K \bsigma(v_h)_{|K} \SCAL \GRAD L_F^K((w_K-w_{\partial K})_{|F}) 
+ (\DIV\bsigma(v_h)_{|K})L_F^K((w_K-w_{\partial K})_{|F})\Big) \dif x \\ 
&=\int_{\partial K} \bsigma(v_h)_{|K}\SCAL\bn_KL_F^K((w_K-w_{\partial K})_{|F})\dif s
= \int_F \bsigma(v_h)_{|K}\SCAL\bn_K(w_K-w_{\partial K})\dif s,
\end{align*}
where we used the divergence formula in $K$. Therefore,
we obtain 
\begin{align*}
n_\sharp(v_h,\hat w_h)) &= \sum_{K\in\calT_h}
\int_{\partial K} \bsigma(v_h)_{|K}\SCAL\bn_K(w_K-w_{\partial K})\dif s \\
&= - \sum_{K\in\calT_h} \lambda_K
\int_{\partial K} \GRAD v_{h|K}\SCAL\bn_K(w_K-w_{\partial K})\dif s \\
&= \sum_{K\in\calT_h} \lambda_K \int_K \big( \GRAD v_{h|K}\SCAL\GRAD(\opRec_K^{k+1}(\hat w_K)-w_K)\big) \dif x,
\end{align*} 
where we used the definition~\eqref{eq:def_Rec_HHO} of the local
reconstruction operator $\opRec_K^{k+1}$ with the test function
$v_{h|K}\in \polP_{k,d}\subset \polP_{k+1,d}$.  
\\
\textup{(ii)} Let us now prove~\eqref{eq:magic_formula2_HHO}. Let $v\in V\loS$
and $\hat w_h\in \hat V_{h,0}^k$. We are going to proceed as in the proof
of~\eqref{eq:n_sharp_2}. We consider the mollification
operators $\calK_\delta\upd:\bL^1(\Dom) \to \bC^\infty(\overline\Dom)$ and
$\calK_\delta\upb:L^1(\Dom) \to C^\infty(\overline\Dom)$
introduced in~\citep[\S3.2]{ErnGu:16_molli}. Let us consider the 
mollified bilinear form
\[ 
n_{\sharp\delta}(v,\hat w_h) := \sum_{K\in\calT_h} \sum_{F\in\calFK} \langle (\calK_\delta\upd(\bsigma(v))\SCAL\bn_K)_{|F},(w_K-w_{\partial K})_{|F}\rangle_F.
\]
By using~\eqref{eq:def_comp_normale} and invoking
the approximation properties of the mollification operators and the
commuting property~\eqref{eq:commut_contrast}, we infer that
$\lim_{\delta\to0} n_{\sharp\delta}(v,\hat w_h)=n_\sharp(v,\hat w_h)$.  
Since the restriction of $\calK_\delta\upd(\bsigma(v))$ to each mesh cell is smooth and since $\calK_\delta\upd(\bsigma(v))\in \bC^0(\overline\Dom)$, we infer that
\begin{align*}
n_{\sharp\delta}(v,\hat w_h) &= \sum_{K\in\calT_h} \int_{\partial K} \calK_\delta\upd(\bsigma(v))\SCAL\bn_K(w_K-w_{\partial K})\dif s = \sum_{K\in\calT_h} \int_{\partial K} \calK_\delta\upd(\bsigma(v))\SCAL\bn_K w_K\dif s \\
&= \sum_{K\in\calT_h} \int_K \big( \calK_\delta\upd(\bsigma(v)) \SCAL \GRAD w_K + \calK_\delta\upb(\DIV\bsigma(v)) w_K\big) \dif x,
\end{align*}
where we used the divergence formula and the commuting
property~\eqref{eq:commut_contrast} in the last line. Letting
$\delta\to0$, we conclude that $n_{\sharp\delta}(v,\hat w_h)$ also
tends to the right-hand side of~\eqref{eq:magic_formula2_HHO} as
$\delta\to0$. Hence, \eqref{eq:magic_formula2_HHO} holds true.\\
\textup{(iii)} The proof of~\eqref{eq:bnd_nsharp_HHO} uses the same
arguments as the proof of Lemma~\ref{lem:bnd_n_sharp}.
\end{proof}

\bRem[\eqref{eq:magic_formula2_HHO}] The right-hand side
of~\eqref{eq:magic_formula2_HHO} does not depend on the face-based
functions $w_{\partial K}$. This identity will replaces the
argument in \citep{DiPEr:15,DiPEL:14} invoking the continuity of the
normal component of $\bsigma(u)$ at the mesh interfaces, which makes
sense only when the exact solution is smooth enough, say
$\bsigma(u)\in \bH^{r}(\Dom)$ with $r>\frac12$.
\end{Rem}

Let $\Vshs := V\loS+P\upb_{k+1}(\calT_h)$ be equipped with the
seminorm $\|v\|_{\Vshs} := |v|_{\lambda,p,q}$ 
defined in~\eqref{eq:snorme_lpq}. Notice that $\|v\|_{\Vshs}=0$
implies that $v=0$ if $v$ has zero mean-value in each mesh cell
$K\in\calT_h$; this is the case for instance if one takes
$v=u-\calE_h^{k+1}(u)$.  We define the consistency error
$\delta_h : \hat V^k_{h,0} \to (\hat V^k_{h,0})'$ by setting,
for all $\hat w_h\in V_{h,0}^k$,
\begin{equation}
\langle \delta_{h}(\hat v_h),\hat w_h\rangle_{(\hat V_{h,0}^k)',\hat V_{h,0}^k} 
:= \hat \ell_h(\hat w_h)-\hat a_h(\hat v_h,\hat w_h).
\end{equation} 
We define global counterparts of the local operators $\opRec_K^{k+1}$, 
$\hat\inter_K^k$, and $\calE_K^{k+1}$, namely $\opRec_h^{k+1}:
\hat V_{h,0}^k \to P\upb_{k+1}(\calT_h)$, 
$\hat\inter_h^k:\Hun\to \hat V^k_{h,0}$, and
$\calE_h^{k+1}:\Hun\to P\upb_{k+1}(\calT_h)$,  by setting 
$\opRec_h^{k+1}(\hat v_h)_{|K}:=\opRec_K^{k+1}(\hat v_K)$,
$\hat\inter_h^k(v)_{|K}:=\hat\inter_K^k(v_{|K})$,
and $\calE_h^{k+1}(v)_{|K}:=\calE_K^{k+1}(v_{|K})$, for all
$\hat v_h\in \hat V_{h,0}^k$, all
$v\in\Hun$, and all $K\in\calT_h$.

\begin{Lem}[Consistency/boundedness] \label{lem:consist_HHO_contrast}
There is $\omega_\sharp$, uniform \wrt 
$\hinH$, $\lambda$, and $u\in V\loS$, but
depending on $p$ and $q$, \sth
\begin{equation}
\|\delta_{h}(\hat\inter_h^k(u))\|_{(\hat V_{h,0}^k)'} \le \omega_\sharp \, \|u-\calE_h^{k+1}(u)\|_{\Vshs}.
\end{equation}
\end{Lem}

\bproof Since $\bsigma(u)=-\lambda\GRAD u$, $\DIV\bsigma(u)=f$, and
$u\in V\loS$, the identity \eqref{eq:magic_formula2_HHO} yields
$\hat \ell_h(\hat w_h) = \sum_{K\in\calT_h} (f,w_K)_{L^2(K)} = 
\sum_{K\in \calT_h} a_K(u,w_K) + n_\sharp(u,\hat w_h)$, where 
$a_K(u,w_K):=\int_K - \bsigma(u)\SCAL\GRAD w_K\dif x$. 
Using the definition of $\hat a_h$ in \eqref{eq:def_hat_a_HHO}, then the identity 
$\opRec_K^{k+1}\circ\hat\inter_K^k = \calE_K^{k+1}$ (see \eqref{eq:ell_proj_HHO}), 
and finally \eqref{eq:magic_formula1_HHO} with $v_h=\calE_h^{k+1}(u)$, we obtain
\begin{align*}
\hat a_h(\hat\inter_h^k(u),\hat w_h) = {}&
\sum_{K\in\calT_h} a_K(\calE_K^{k+1}(u),w_K) + n_\sharp(\calE_h^{k+1}(u),\hat w_h) \\
&+\sum_{K\in\calT_h} 
\lambda_K (h_{\partial K}^{-1}\opSt_{\partial K}^k(\hat\inter_K^k(u)),
\opSt_{\partial K}^k(\hat w_K))_{L^2(\partial K)}.
\end{align*}
Subtracting these two identities and using the definition of 
$\calE_K^{k+1}(u)$, which implies that $a_K(u-\calE_K^{k+1}(u),w_K)=0$, for 
all $K\in\calT_h$, leads to $\langle\delta_{h}(\hat\inter_h^k(u)),\hat w_h\rangle_{(\hat
V_{h,0}^k)',\hat V_{h,0}^k} = \term_1+\term_2$ with
\begin{align*}
\term_1 &:= n_\sharp(u-\calE_h^{k+1}(u),\hat w_h),\qquad
\term_2 := -
\sum_{K\in\calT_h} \lambda_K (h_{\partial K}^{-1}
\opSt_{\partial K}^k(\hat\inter_K^k(u)),
\opSt_{\partial K}^k(\hat w_K)_{L^2(\partial K)}.
\end{align*}
We invoke~\eqref{eq:bnd_nsharp_HHO} to bound $\term_1$ and observe 
that $\sum_{K\in\calT_h} \lambda_Kh_K^{-1}\|w_K-w_{\partial K}\|_{L^2(\partial K)}^2
\le \|\hat w_h\|_{\hat V_{h,0}^k}^2$ owing to~\eqref{eq:norme_HHO_contraste}.
For the bound on $\term_2$, we proceed as in \citep{DiPEr:15,DiPEL:14}.
\end{proof}

\begin{Th}[Error estimate] \label{th:cv_contrast_HHO}
Let $u$ solve~\eqref{model_elliptic_contrast} and $\hat u_h$
solve~\eqref{eq:weak_HHO_contrast} with $\hat a_h$ and $\hat \ell_h$ defined
in~\eqref{eq:def_hat_a_HHO}. 
Assume that there is $r>0$ \sth $u\in H^{1+r}(\Dom)$. 
There is $c$, uniform \wrt $\hinH$, $\lambda$, and
$u\in H^{1+r}(\Dom)$, 
but depending on $r$, \sth the following holds true:
\begin{equation} \label{eq:err_est_HHO_contrast}
\|\lambda^{\frac12} \GRAD_h (u-\opRec_h^{k+1}(\hat u_h))\|_{\bL^2(\Dom)}
\le c\, \|u-\calE_h^{k+1}(u)\|_{\Vshs}.
\end{equation}
Moreover, letting $t:=\min(r,k+1)$, $\chi_t=1$ if $t\le 1$ and $\chi_t=0$ if
$t>1$, we have
\begin{multline}
\|\lambda^{\frac12} \GRAD_h (u-\opRec_h^{k+1}(\hat u_h))\|_{\bL^2(\Dom)}
\\ \le c\, \bigg( \sum_{K\in\calT_h} \!\!\lambda_Kh_K^{2t}|u|_{H^{1+t}(K)}^2 + 
\frac{\chi_t}{\lambda_K} 
h_K^{2d(\frac{d+2}{2d}-\frac{1}{q})}\|f\|_{L^q(K)}^2\bigg)^{\frac12}\!. 
\label{eq:cv_rate_HHO_contrast}
\end{multline}
\end{Th}
\begin{proof}
\textup{(i)} We adapt the proof of Lemma~\ref{lem:abst_err_est_QO} to 
exploit the convergence order of the reconstruction operator.
Let us set $\hat \zeta_h^k := \hat\inter_h^k(u)-\hat u_h\in \hat V_{h,0}^k$
so that $\hat \zeta_K^k = \hat\inter_K^k(u_{|K})-\hat u_K$ for all $K\in\calT_h$. 
The coercivity property from Lemma~\ref{lem:stab_HHO} 
and the definition of the consistency error imply that 
\begin{align*}
&\alpha\, \|\lambda^{\frac12}\GRAD_h \opRec_h^{k+1}(\hat \zeta_h^k )\|_{\bL^2(\Dom)}^2
\le \frac{\hat a_h(\hat \zeta_h^k,\hat \zeta_h^k)}{\|\hat \zeta_h^k\|_{\hat V_{h,0}^k}^2}
\|\lambda^{\frac12}\GRAD_h \opRec_h^{k+1}(\hat \zeta_h^k )\|_{\bL^2(\Dom)}^2\\
&\le \frac{\big(\hat a_h(\hat \zeta_h^k,\hat \zeta_h^k)\big)^2}{\|\hat \zeta_h^k\|_{\hat V_{h,0}^k}^2}
=\frac{\langle \delta_h(\hat\inter^k_h(u)),\hat \zeta_h^k\rangle^2_{(\hat V_{h,0}^k)',\hat V_{h,0}^k}}{\|\hat \zeta_h^k\|_{\hat V_{h,0}^k}^2}\le \|\delta_h(\hat\inter^k_h(u))\|^2_{(\hat V_{h,0}^k)'}.
\end{align*}
Then, lemma~\ref{lem:consist_HHO_contrast} yields
$\|\lambda^{\frac12}\GRAD \opRec_h^{k+1}(\hat \zeta_h^k
)\|_{\bL^2(\Dom)} \le c\|u-\calE_h^{k+1}(u)\|_{\Vshs}$.  Moreover, since
$\opRec_K^{k+1}(\hat\inter_K^k(u))=\calE_K^{k+1}(u)$ for all $K\in\calT_h$, see
\eqref{eq:ell_proj_HHO}, we have
\[
u-\opRec_h^{k+1}(\hat u_{h})
= u-\calE_h^{k+1}(u) + \opRec_h^{k+1}(\hat \zeta_h^k).
\]
The estimate \eqref{eq:err_est_HHO_contrast} is now a consequence of the triangle inequality.\\
\textup{(ii)} We now prove \eqref{eq:cv_rate_HHO_contrast}. 
Let us set $\eta^{k+1}:=u-\calE_h^{k+1}(u)$. We need to bound
$\|\eta^{k+1}\|_{\Vshs} = |\eta^{k+1}|_{\lambda,p,q}$, \ie we must estimate 
$\|\GRAD\eta^{k+1}\|_{\bL^2(K)}$,
$h_K^{d(\frac12-\frac1p)}\|\GRAD\eta^{k+1}\|_{\bL^p(K)}$, and
$h_K^{d(\frac{d+2}{2d}-\frac{1}{q})}\|\Delta\eta^{k+1}\|_{L^q(K)}$
(see~\eqref{eq:snorme_lpq}).
Owing to the optimality property of the elliptic projection and the
approximation properties of $\Pi_K^{k+1}$, we have
\[
\|\GRAD\eta^{k+1}\|_{\bL^2(K)} \le \|\GRAD(u-\Pi_K^{k+1}(u))\|_{\bL^2(K)} \le c\, h_K^{t}|u|_{H^{1+t}(K)}.
\]
for $t=\min(r,k+1)$.
Let us now consider the other two terms. 
Let $\ell:=\lceil t\rceil$, so that 
 $t\le \ell\le 1+t$. Notice also that $\ell\le k+1$, and 
$\ell\ge 1$ since we assumed that $r>0$. Let us set
$\eta^\ell:=u-\calE_h^\ell(u)$, then
$\|\GRAD\eta^{\ell}\|_{\bL^2(K)}\le ch_K^t|u|_{H^{1+t}(K)}$.
Invoking the triangle inequality,
an inverse inequality, and the triangle inequality again, we infer that
\begin{align*}
h_K^{d(\frac12-\frac1p)}\|\GRAD\eta^{k+1}\|_{\bL^p(K)}
\le {}&h_K^{d(\frac12-\frac1p)}\|\GRAD\eta^{\ell}\|_{\bL^p(K)}
+ c\, \big(\|\GRAD\eta^{k+1}\|_{\bL^2(K)}+\|\GRAD\eta^{\ell}\|_{\bL^2(K)}\big),
\end{align*}
and the two terms between the parentheses are bounded by $ch_K^t|u|_{H^{1+t}(K)}$.
Moreover, invoking~\eqref{eq:inv_ineq_Lp_Hr_contrast}, we obtain
\begin{align*}
h_K^{d(\frac12-\frac{1}{p})}\|\GRAD\eta^\ell\|_{\bL^p(K)} &\le c\,
\big(\|\GRAD\eta^\ell\|_{\bL^2(K)} + h_{K}^t
|\GRAD\eta^\ell|_{\bH^t(K)}\big) \\
&=c\,\big(\|\GRAD\eta^\ell\|_{\bL^2(K)} + h_{K}^t
|u|_{H^{1+t}(K)}\big) \le c'\,h_K^t|u|_{H^{1+t}(K)},
\end{align*}
since $t\le \ell$. Similarly, we have
\begin{align*} 
h_K^{d(\frac{d+2}{2d}-\frac{1}{q})}\|\Delta\eta^{k+1}\|_{L^q(K)}
\le {}& h_K^{d(\frac{d+2}{2d}-\frac{1}{q})}\|\Delta\eta^{\ell}\|_{L^q(K)}
+ c\, \big(\|\GRAD\eta^{k+1}\|_{\bL^2(K)}+\|\GRAD\eta^{\ell}\|_{\bL^2(K)}\big).
\end{align*}
It remains to estimate 
$h_K^{d(\frac{d+2}{2d}-\frac{1}{q})}\|\Delta\eta^{\ell}\|_{L^q(K)}$. We proceed
as in the end of the proof of Theorem~\ref{th:conv_Nitsche_contrasted}. 
If $t\le 1$ (so that $\chi_t=1$), we have $\ell=1$, and we infer that
\[
h_K^{d(\frac{d+2}{2d}-\frac{1}{q})}\|\Delta\eta^{\ell}\|_{L^q(K)}
= \lambda_K^{-1}h_K^{d(\frac{d+2}{2d}-\frac{1}{q})}\|f\|_{L^q(K)}.
\]
Otherwise, we have $t>1$ (so that $\chi_t=0$) and $\ell\ge2$, and we 
take $q=2$. Then,
using the triangle inequality, an inverse inequality, 
and the triangle inequality again, we obtain
\begin{align*}
h_K\|\Delta\eta^{\ell}\|_{L^q(K)}
\le {}&h_K\|\Delta(u-\Pi_K^\ell(u)\|_{L^q(K)}
\\ &+ c\,\big(\|\GRAD(u-\Pi_K^\ell(u))\|_{\bL^2(K)}+\|\GRAD\eta^{\ell}\|_{\bL^2(K)}\big),
\end{align*}
where $\Pi_K^\ell$ is the $L^2$-orthogonal projection 
onto $\polP_{\ell,d}$. We conclude by invoking the approximation 
properties of $\Pi_K^\ell$, recalling that 
$\|\GRAD\eta^{\ell}\|_{\bL^2(K)}\le ch_K^t|u|_{H^{1+t}(K)}$.
\end{proof}

\bRem[Supercloseness]
Step~(i) in the above proof actually shows that 
$\|\hat \zeta_h^k\|_{\hat V_{h,0}^k}\le c\|u-\calE_h^{k+1}(u)\|_{\Vshs}$.
Since $\zeta_K^k=\Pi_K^k(u)-u_K$ for all $K\in\calT_h$, 
this implies the supercloseness bound
$(\sum_{K\in\calT_h} \lambda_K\|\GRAD(\Pi_K^k(u)-u_K)\|_{\bL^2(K)}^2)^{\frac12}
\le c\|u-\calE_h^{k+1}(u)\|_{\Vshs}$. 
\eRem

\section{Extensions to Maxwell's equations}
\label{sec:Maxwell}

The various techniques presented in this paper can be extended to the
context of Maxwell's equations, since arguments similar to those
exposed in \S\ref{Sec:nsharp} can be deployed to define the tangential
trace of vectors fields on a face of $K$. Without going into the details,
we show in this section how that can be done.

\subsection{Lifting and tangential trace}
Let  $p$, $q$ be real numbers satisfying
\eqref{eq:p_q_face_to_cell}, and let $\tp\in (2,p]$ be such that
$q\ge \frac{\tp d}{\tp +d}$.  Let $K$ be a cell in $\calT_h$, and
let $F\in\calF_K$ be a face of $K$.  Following
\citep{Ern_Guermond_cmam_2017}, we introduce the space
\begin{align}
\bY\upc(F) := \{\bphi\in \bW^{\frac{1}{\tp},{\tp'}}(F)\tq 
\bphi\SCAL\bn_F=0\}, \label{eq:def_Yupc_F}
\end{align}
which we equip with the norm
$\|\bphi\|_{\bY\upc(F)}:= \|\bphi\|_{\bL^{\tp'}(F)}+
h_F^{\frac{1}{\tp}} |\bphi|_{\bW^{\frac{1}{\tp},\tp'}(F)}$.  Then the
following result can be established by proceeding as in the proof of
Lemma~\ref{lem:lifting_face_to_cell}.

\begin{Lem}[Face-to-cell Lifting]
There exist a constant $c$, 
uniform \wrt $h$, but depending on $p$ and $q$, 
and a lifting
operator $E_F^K : \bY\upc(F) \rightarrow \bW^{1,\tp'}(K)$ such
that the following holds true for any $\bphi\in \bY\upc(F)$:
$E_F^K(\bphi)_{|\partial K{\setminus} F}=\bzero$,
$E_F^K(\bphi)_{|F}=\bphi$, and
\begin{equation} \label{eq:stab_lifting_Vc_p_q}
|E_F^K(\bphi)|_{\bW^{1,p'}(K)} 
+ h_K^{-1+d(\frac1q-\frac1p)}\|E_F^K(\bphi)\|_{\bL^{q'}(K)} 
\le c\, h_K^{-\frac{1}{\tp} + d(\frac{1}{\tp}-\frac{1}{p})} \|\bphi\|_{\bY\upc(F)}.
\end{equation}
\end{Lem}

With this lifting operator in hand, we can 
define an extension to the notion
of the tangential trace on $F$ of a vector field. To this end,
we introduce the functional space
\begin{align}
\bS\upc(K) := \{ \btau\in\bL^p(K) \tq \ROT\btau \in \bL^q(K) \},
\end{align}
where the superscript $\upc$ refers to the fact that the tangential trace is related to
the curl operator. We equip $\bS\upc(K)$ with the following dimensionally-consistent norm:
\begin{equation}
\|\btau\|_{\bS\upc(K)}:=\|\btau\|_{\bL^p(K)} +
  h_{K}^{1+d(\frac1p-\frac1q)}\|\ROT\btau\|_{L^q(K)}. \label{def_of_bSupC_K}
\end{equation}
We now define the tangential trace of any field
$\btau$ in $\bS\upc(K)$ on the face $F$ of $K$ 
to be the linear form
$(\btau\CROSS \bn_K)_{|F}\in \bY\upc(F)'$ such that
\begin{equation} \label{eq:def_dual_tang_comp}
\langle (\btau\CROSS \bn_K)_{|F},\bphi \rangle_F :=
\int_K \bigg( \btau\SCAL\ROT E_F^K(\bphi) - (\ROT\btau)\SCAL E_F^K(\bphi)\bigg) \dif x,
\end{equation}
for all $\bphi\in \bY\upc(F)$, where $\langle\SCAL,\SCAL\rangle_F$ 
now denotes the duality pairing between $\bY\upc(F)'$ and $\bY\upc(F)$. 
Note that the right-hand side of~\eqref{eq:def_dual_tang_comp} 
is well-defined owing to H\"older's inequality and~\eqref{eq:stab_lifting_Vc_p_q}.

The discretization now involves the vector-valued broken finite element space
\begin{equation}
\bP\upb_k(\calT_h) = \bset \bv_h\in
\bL^\infty(\Dom) \tq \bv_{h|K}\in \bP_K,\, \forall K\in\calT_h\eset,
\end{equation}
where $\bP_K:=(\mapK)^{-1}(\wbP)\subset \bW^{k+1,\infty}(K)$, 
$\wbKPS$ is the reference element, and $\mapK$
is an appropriate transformation. For instance, one can take
$\mapK(\bv)=\mapKg(\bv):=\bv\circ \bT_K$ for continuous Lagrange
elements or for dG approximation; one can also take
$\mapK(\bv)=\mapKc(\bv):=\Jac_K\tr (\bv\circ\bT_K)$ for edge elements
($\mapKc$ is covariant Piola transformation and $\Jac_K$ the 
Jacobian of the geometric mapping). For any face
$F\in \calF_K$, we denote by $\bP_F$ the trace of $\bP_K$ on $F$. The
following result is the counterpart of
Lemma~\ref{lem:bnd_normal_contrast}.

\begin{Lem}[Bound on tangential component]\label{lem:bnd_tang_contrast}
There exists a constant $c$, uniform \wrt $h$,
but depending on $p$ and $q$,
so that the following 
estimate holds true for all $\bv\in \bS\upc(K)$, 
\begin{equation} \label{eq:bnd_tang_dual}
  \|(\bv\CROSS\bn_K)_{|F}\|_{Y\upc(F)'} \le c\,
  h_K^{-\frac{1}{\tp}+d(\frac{1}{\tp}-\frac{1}{p})} \|\bv\|_{\bS\upc(K)}.
\end{equation}
Moreover, we have
\begin{equation} \label{eq:bnd_tang_dual_h} |\langle
  (\bv\CROSS\bn_K)_{|F},\bphi_h \rangle| \le c\,
  h_{K}^{d(\frac12-\frac1p)}\|\bv\|_{\bS\upc(K)}
  h_F^{-\frac12}\|\bphi_h\|_{L^2(F)},
\end{equation}
for all $\bphi_h\in \bP_F$ s.t.~$\bphi\SCAL\bn_F=0$, all $K\in\calT_h$, and all $F\in\calF_K$.
\end{Lem}

Lemma~\ref{lem:bnd_tang_contrast} is essential 
for the error analysis of nonconforming
approximation techniques of Maxwell's equations.  It is a
generalization of \cite[Lem.~A3]{Bonito_Guermond_Luddens_M2aN_2016}
and \cite[Lem.~8.2]{Buffa_Perugia_SINUM_2006}.

\subsection{Definition of $n_\sharp\upc$ and key identities}

The consistency analysis of Nitsche's boundary penalty method and of the
dG approximation applied to Maxwel's equations can be done by introducing a
bilinear form $n_\sharp$ as in \S\ref{Sec:nsharp}. We henceforth assume
that the space dimension is either $d=2$ or $d=3$.

We define the notion of diffusive flux by introducing
$\bsigma: \Hrot\to \bL^2(\Dom)$ such that
$\bsigma(\bv):=\lambda \ROT \bv$, for any
$\bv\in\Hrot$. Here, the diffusivity $\lambda$ is either the
reciprocal of the magnetic permeability or the reciprocal of
electrical conductivity, depending whether one works with the electric
field or the magnetic field. The diffusivity is assumed to satisfy
the hypotheses introduced in Section~\ref{Sec:Preliminaries}.
We further define 
\begin{equation}
\bV\loS := \{ \bv\in \Hrot\st \bsigma(\bv)\in \bL^p(\Dom), \ \ROT\bsigma(\bv)\in \bL^q(\Dom)\},
\end{equation}
and set $\bV_\sharp := \bV\loS+\bP\upb_k(\calT_h)$.

We adopt the same notation as in \S\ref{Sec:nsharp}.
Recall that for  any $K\in \calT_h$ and any $F\in\calF_K$, we have defined
$\epsilon_{K,F}=\bn_F\SCAL\bn_K=\pm1$. We consider arbitrary weights 
$\theta_{K,F}$ satisfying \eqref{convex_interface_weights}.
We introduce the bilinear form 
$n_\sharp\upc:(\bV\loS+\bP\upb_k(\calT_h)) \times \bP\upb_k(\calT_h)\to \Real$ defined as
follows:
\begin{align} \label{eq:def_n_sharp_ROT}
n_\sharp\upc(\bv,\bw_h) := {}&\sum_{F\in\calFh} \sum_{K\in\calT_F}
\epsilon_{K,F}\theta_{K,F} 
\langle (\bsigma(\bv)_{|K}\CROSS\bn_K)_{|F},\jump{\Pi_F(\bw_h)} \rangle_F,
\end{align}
where $\Pi_F$ is the $\ell^2$-orthogonal projection onto the
hyperplane tangent to $F$, \ie
$\Pi_F(\bb_h):=\bb_h - (\bb_h\SCAL \bn_K) \bn_K =\bn_K\CROSS(\bb_h\CROSS \bn_K)$. Notice that~\eqref{eq:def_n_sharp_ROT} is
  meaningful since $\Pi_F(\bb_h)_{|F}$ is in
  $\bW^{\frac{1}{\tp},{\tp'}}(F)$ and
$\Pi_F(\bb_h)\SCAL\bn_F=0$, \ie $\Pi_F(\bb_h)\in \bY\upc(F)$ for any $F\in\calF_h$. The following result is the counterpart of Lemma~\ref{lem:identity_n_sharp_dif}.

\begin{Lem}[Identities for $n_\sharp\upc$] 
\label{lem:identity_n_sharp_dif_ROT}%
The following holds true for any choice of weights
$\{\theta_{K,F}\}_{F\in \calFh,K\in\calT_F}$ and for all
$\bw_h\in \bP\upb_k(\calT_h)$, all $\bv_h\in \bP\upb_k(\calT_h)$, and all $\bv\in \bV\loS$:
\begin{subequations}\begin{align} \label{eq:n_sharp_1_ROT}%
n_\sharp\upc(\bv_h,\bw_h) 
&= \sum_{F\in\calFh} \int_F (\avg{\bsigma(\bv_h)}_\theta\CROSS\bn_F) \SCAL\jump{\Pi_F(\bw_h)} \dif s,\\
\label{eq:n_sharp_2_ROT}
n_\sharp\upc(\bv,\bw_h) &= \sum_{K\in\calT_h} 
\int_K \Big(\bsigma(v)\SCAL\ROT \bw_{h|K}
- (\ROT\bsigma(v)) \SCAL \bw_{h|K}\Big) \dif x.
\end{align}\end{subequations}
\end{Lem}
\begin{proof} 
The proof is similar to that of
Lemma~\ref{lem:identity_n_sharp_dif}.  The proof of
\eqref{eq:n_sharp_1_ROT} is quasi-identical to that of
\eqref{eq:n_sharp_1}.  For the proof of \eqref{eq:n_sharp_2_ROT}, one
invokes the mollifying operators
$\calK_\delta\upc:\bL^1(\Dom) \to \bC^\infty(\overline\Dom)$ and
$\calK_\delta\upd:\bL^1(\Dom) \to \bC^\infty(\overline\Dom)$
introduced in~\citep[\S3.2]{ErnGu:16_molli}. 
These two operators satisfy the following
key commuting property:
\begin{equation} \label{eq:commut_contrast_ROT}
\ROT(\calK_\delta\upc(\btau)) = \calK_\delta\upd(\ROT\btau),
\end{equation}
for all $\btau\in \bL^1(\Dom)$ s.t.~$\ROT\btau\in \bL^1(\Dom)$.  Then
one uses the identities
$\jump{\bv\CROSS\Pi_F(\bw)} = \avg{\bv}_{\theta}\CROSS
\jump{\Pi_F(\bw)} + \jump{\bv}\CROSS\avg{\Pi_F(\bw)}_{\bar\theta}$,
$\bn_K\CROSS\Pi_F(\bw_h)=\bn_K\CROSS \bw_h$, and
$\DIV(\bw_h\CROSS \bsigma(\bv))=\bsigma(\bv)\SCAL(\ROT\bw_h) -
\bw_h\SCAL (\ROT \bsigma(\bv))$.
\end{proof}

We now establish the boundedness of the bilinear
form $n_\sharp\upc$. Since $\bsigma(\bv)_{|K}\in \bS\upc(K)$ 
for all $K\in\calT_h$ 
and all $\bv\in \bV\loS+\bP\upb_k(\calT_h)$, we equip the space 
$\bV\loS+\bP\upb_k(\calT_h)$ with the seminorm
\begin{multline}
|\bv|_{n_\sharp\upc}^2 := \sum_{K\in \calT_h} 
\lambda_K^{-1}\Big(h_K^{2d(\frac12-\frac1p)} \|\bsigma(\bv)_{|K}\|_{\bL^p(K)}^2 \\
+h_K^{2d(\frac{2+d}{2d}-\frac{1}{q})} 
\|\ROT\bsigma(\bv)_{|K}\|_{\bL^q(K)}^2\Big).\label{eq:seminorm_sharp_n_ROT}
\end{multline}%

\begin{Lem}[Boundedness of $n_\sharp\upc$] \label{lem:bnd_n_sharp_ROT}%
With the weights defined in \eqref{eq:def_weight_theta} and 
$\lambda_F$ defined in \eqref{def_of_lambdaF} for all $F\in\calFh$,
there is $c$, uniform \wrt $\hinH$ and $\lambda$, but 
depending on $p$ and $q$, \sth
the following holds true for all
$\bv\in \bV\loS+\bP\upb_k(\calT_h)$ and all $\bw_h\in \bP\upb_k(\calT_h)$:
\begin{equation}
|n_\sharp\upc(\bv,\bw_h)| \le c\, |\bv|_{n_\sharp\upc} \bigg( \sum_{F\in\calFh} 
\lambda_F h_F^{-1} \|\jump{\Pi_F(\bw_h)}\|_{\bL^2(F)}^2 \bigg)^{\frac12}. \label{eq1:lem:bnd_n_sharp_ROT}
\end{equation}
\end{Lem}

With the above tools in hand, one can revisit
\cite{Buffa_Perugia_SINUM_2006} and greatly simplify the analysis of
the dG approximation of Maxwell's equations.  One can also extend the
work in \citep{Ern_Guermond_CAMWA_2018} and analyze Nitsche's boundary
penalty technique with edge elements; one can also revisit
\cite{Bonito_guermond_Luddens_amd2_2013}, where Nitsche's boundary penalty
technique has been used in conjunction with Lagrange elements.  In all
the cases one then obtains error estimates that are robust with
respect to the diffusivity contrast.

\bibliographystyle{abbrvnat}

\bibliography{ref}

\begin{thebibliography}{35}
\providecommand{\natexlab}[1]{#1}
\providecommand{\url}[1]{\texttt{#1}}
\expandafter\ifx\csname urlstyle\endcsname\relax
  \providecommand{\doi}[1]{doi: #1}\else
  \providecommand{\doi}{doi: \begingroup \urlstyle{rm}\Url}\fi

\bibitem[Amrouche et~al.(1998)Amrouche, Bernardi, Dauge, and Girault]{AmBDG:98}
C.~Amrouche, C.~Bernardi, M.~Dauge, and V.~Girault.
\newblock Vector potentials in three-dimensional non-smooth domains.
\newblock \emph{Math. Methods Appl. Sci.}, 21\penalty0 (9):\penalty0 823--864,
  1998.

\bibitem[Arnold(1982)]{Arnol:82}
D.~N. Arnold.
\newblock An interior penalty finite element method with discontinuous
  elements.
\newblock \emph{SIAM J. Numer. Anal.}, 19:\penalty0 742--760, 1982.

\bibitem[Badia et~al.(2014)Badia, Codina, Gudi, and Guzm\'an]{BaCGG:14}
S.~Badia, R.~Codina, T.~Gudi, and J.~Guzm\'an.
\newblock Error analysis of discontinuous {G}alerkin methods for the {S}tokes
  problem under minimal regularity.
\newblock \emph{IMA J. Numer. Anal.}, 34\penalty0 (2):\penalty0 800--819, 2014.

\bibitem[Bernardi and Girault(1998)]{Bernardi_Girault_1998}
C.~Bernardi and V.~Girault.
\newblock A local regularization operator for triangular and quadrilateral
  finite elements.
\newblock \emph{SIAM J. Numer. Anal.}, 35\penalty0 (5):\penalty0 1893--1916,
  1998.

\bibitem[Bernardi and Hecht(2002)]{BerHe:02}
C.~Bernardi and F.~Hecht.
\newblock Error indicators for the mortar finite element discretization of the
  {L}aplace equation.
\newblock \emph{Math. Comp.}, 71\penalty0 (240):\penalty0 1371--1403, 2002.

\bibitem[Bernardi and Verf\"urth(2000)]{BerVe:00}
C.~Bernardi and R.~Verf\"urth.
\newblock Adaptive finite element methods for elliptic equations with
  non-smooth coefficients.
\newblock \emph{Numer. Math.}, 85\penalty0 (4):\penalty0 579--608, 2000.

\bibitem[Bonito et~al.(2013)Bonito, Guermond, and
  Luddens]{Bonito_guermond_Luddens_amd2_2013}
A.~Bonito, J.-L. Guermond, and F.~Luddens.
\newblock Regularity of the maxwell equations in heterogeneous media and
  lipschitz domains.
\newblock \emph{J. Math. Anal. Appl.}, 408:\penalty0 498--512, 2013.

\bibitem[Bonito et~al.(2016)Bonito, Guermond, and
  Luddens]{Bonito_Guermond_Luddens_M2aN_2016}
A.~Bonito, J.-L. Guermond, and F.~Luddens.
\newblock An interior penalty method with {$C^0$} finite elements for the
  approximation of the {M}axwell equations in heterogeneous media: convergence
  analysis with minimal regularity.
\newblock \emph{ESAIM Math. Model. Numer. Anal.}, 50\penalty0 (5):\penalty0
  1457--1489, 2016.

\bibitem[Buffa and Perugia(2006)]{Buffa_Perugia_SINUM_2006}
A.~Buffa and I.~Perugia.
\newblock Discontinuous {G}alerkin approximation of the {M}axwell eigenproblem.
\newblock \emph{SIAM J. Numer. Anal.}, 44\penalty0 (5):\penalty0 2198--2226,
  2006.

\bibitem[Burman and Zunino(2006)]{BurZu:06}
E.~Burman and P.~Zunino.
\newblock A domain decomposition method for partial differential equations with
  non-negative form based on interior penalties.
\newblock \emph{SIAM J. Numer. Anal.}, 44:\penalty0 1612--1638, 2006.

\bibitem[Cai et~al.(2011)Cai, Ye, and Zhang]{CaiHZ:11}
Z.~Cai, X.~Ye, and S.~Zhang.
\newblock Discontinuous {G}alerkin finite element methods for interface
  problems: a priori and a posteriori error estimations.
\newblock \emph{SIAM J. Numer. Anal.}, 49\penalty0 (5):\penalty0 1761--1787,
  2011.

\bibitem[Carstensen and Schedensack(2015)]{CarSc:15}
C.~Carstensen and M.~Schedensack.
\newblock Medius analysis and comparison results for first-order finite element
  methods in linear elasticity.
\newblock \emph{IMA J. Numer. Anal.}, 35\penalty0 (4):\penalty0 1591--1621,
  2015.

\bibitem[Cockburn et~al.(2016)Cockburn, Di~Pietro, and Ern]{CoDPE:16}
B.~Cockburn, D.~A. Di~Pietro, and A.~Ern.
\newblock Bridging the {Hybrid High-Order} and {Hybridizable Discontinuous
  Galerkin} methods.
\newblock \emph{ESAIM: Math. Model Numer. Anal. (M2AN)}, 50\penalty0
  (3):\penalty0 635--650, 2016.

\bibitem[Crouzeix and Raviart(1973)]{CR73}
M.~Crouzeix and P.-A. Raviart.
\newblock Conforming and nonconforming finite element methods for solving the
  stationary {S}tokes equations. {I}.
\newblock \emph{Rev. Fran\c caise Automat. Informat. Recherche Op\'erationnelle
  S\'er. Rouge}, 7\penalty0 (R-3):\penalty0 33--75, 1973.

\bibitem[Di~Pietro and Ern(2012)]{DiPEr:12}
D.~A. Di~Pietro and A.~Ern.
\newblock \emph{Mathematical {A}spects of {D}iscontinuous {G}alerkin
  {M}ethods}, volume~69 of \emph{Math\'ematiques \& Applications}.
\newblock Springer-Verlag, Berlin, 2012.

\bibitem[Di~Pietro and Ern(2015)]{DiPEr:15}
D.~A. Di~Pietro and A.~Ern.
\newblock A {Hybrid High-Order} locking-free method for linear elasticity on
  general meshes.
\newblock \emph{Comput. Meth. Appl. Mech. Engrg.}, 283:\penalty0 1--21, 2015.

\bibitem[Di~Pietro et~al.(2008)Di~Pietro, Ern, and Guermond]{DiErG:08}
D.~A. Di~Pietro, A.~Ern, and J.-L. Guermond.
\newblock Discontinuous {G}alerkin methods for anisotropic semi-definite
  diffusion with advection.
\newblock \emph{SIAM J. Numer. Anal.}, 46\penalty0 (2):\penalty0 805--831,
  2008.

\bibitem[Di~Pietro et~al.(2014)Di~Pietro, Ern, and Lemaire]{DiPEL:14}
D.~A. Di~Pietro, A.~Ern, and S.~Lemaire.
\newblock An arbitrary-order and compact-stencil discretization of diffusion on
  general meshes based on local reconstruction operators.
\newblock \emph{Comput. Meth. Appl. Math.}, 14\penalty0 (4):\penalty0 461--472,
  2014.

\bibitem[Dryja(2003)]{Dryja:03}
M.~Dryja.
\newblock On discontinuous {G}alerkin methods for elliptic problems with
  discontinuous coefficients.
\newblock \emph{Comput. Methods Appl. Math.}, 3\penalty0 (1):\penalty0 76--85,
  2003.

\bibitem[Dryja et~al.(2007)Dryja, Galvis, and Sarkis]{DryGS:07}
M.~Dryja, J.~Galvis, and M.~Sarkis.
\newblock B{DDC} methods for discontinuous {G}alerkin discretization of
  elliptic problems.
\newblock \emph{J. Complexity}, 23\penalty0 (4-6):\penalty0 715--739, 2007.

\bibitem[Ern and Guermond(2006)]{ErnGu:06}
A.~Ern and J.-L. Guermond.
\newblock Discontinuous {G}alerkin methods for {F}riedrichs' systems. {I}.
  {G}eneral theory.
\newblock \emph{SIAM J. Numer. Anal.}, 44\penalty0 (2):\penalty0 753--778,
  2006.

\bibitem[Ern and Guermond(2016)]{ErnGu:16_molli}
A.~Ern and J.-L. Guermond.
\newblock Mollification in strongly {L}ipschitz domains with application to
  continuous and discrete de {R}ham complexes.
\newblock \emph{Comput. Methods Appl. Math.}, 16\penalty0 (1):\penalty0 51--75,
  2016.

\bibitem[Ern and Guermond(2017)]{Ern_Guermond_M2AN_2017}
A.~Ern and J.-L. Guermond.
\newblock Finite element quasi-interpolation and best approximation.
\newblock \emph{M2AN Math. Model. Numer. Anal.}, 51\penalty0 (4):\penalty0
  1367--1385, 2017.

\bibitem[Ern and Guermond(2018{\natexlab{a}})]{Ern_Guermond_CAMWA_2018}
A.~Ern and J.-L. Guermond.
\newblock Analysis of the edge finite element approximation of the {M}axwell
  equations with low regularity solutions.
\newblock \emph{Comput. Math. Appl.}, 75\penalty0 (3):\penalty0 918--932,
  2018{\natexlab{a}}.

\bibitem[Ern and Guermond(2018{\natexlab{b}})]{Ern_Guermond_cmam_2017}
A.~Ern and J.-L. Guermond.
\newblock Abstract nonconforming error estimates and application to boundary
  penalty methods for diffusion equations and time-harmonic {M}axwell's
  equations.
\newblock \emph{Comput. Methods Appl. Math.}, 18\penalty0 (3):\penalty0
  451--475, 2018{\natexlab{b}}.

\bibitem[Ern et~al.(2009)Ern, Stephansen, and Zunino]{ErStZ:09}
A.~Ern, A.~F. Stephansen, and P.~Zunino.
\newblock A discontinuous {G}alerkin method with weighted averages for
  advection-diffusion equations with locally small and anisotropic diffusivity.
\newblock \emph{IMA J. Numer. Anal.}, 29\penalty0 (2):\penalty0 235--256, 2009.

\bibitem[Gagliardo(1957)]{Gagliardo:57}
E.~Gagliardo.
\newblock Caratterizzazioni delle tracce sulla frontiera relative ad alcune
  classi di funzioni in {$n$} variabili.
\newblock \emph{Rend. Sem. Mat. Univ. Padova}, 27:\penalty0 284--305, 1957.

\bibitem[Grisvard(1985)]{Gr85}
P.~Grisvard.
\newblock \emph{Elliptic problems in nonsmooth domains}, volume~24 of
  \emph{Monographs and Studies in Mathematics}.
\newblock Pitman (Advanced Publishing Program), Boston, MA, 1985.

\bibitem[Gudi(2010)]{Gudi:10}
T.~Gudi.
\newblock A new error analysis for discontinuous finite element methods for
  linear elliptic problems.
\newblock \emph{Math. Comp.}, 79\penalty0 (272):\penalty0 2169--2189, 2010.

\bibitem[Jochmann(1999)]{Jochmann_1999}
F.~Jochmann.
\newblock An ${H}^s$-regularity result for the gradient of solutions to
  elliptic equations with mixed boundary conditions.
\newblock \emph{J. Math. Anal. Appl.}, 238:\penalty0 429--450, 1999.

\bibitem[Li and Mao(2013)]{Li_Shipeng_2013}
M.~Li and S.~Mao.
\newblock A new a priori error analysis of nonconforming and mixed finite
  element methods.
\newblock \emph{Appl. Math. Lett.}, 26\penalty0 (1):\penalty0 32--37, 2013.

\bibitem[Nitsche(1971)]{Nitsche_1971}
J.~Nitsche.
\newblock \"{U}ber ein {V}ariationsprinzip zur {L}\"osung von
  {D}irichlet-{P}roblemen bei {V}erwendung von {T}eilr\"aumen, die keinen
  {R}andbedingungen unterworfen sind.
\newblock \emph{Abh. Math. Sem. Univ. Hamburg}, 36:\penalty0 9--15, 1971.

\bibitem[Sch{\"o}berl(2001)]{Schoberl_2001}
J.~Sch{\"o}berl.
\newblock Commuting quasi-interpolation operators for mixed finite elements.
\newblock Technical Report ISC-01-10-MATH, Texas A\&M University, 2001.
\newblock URL \url{www.isc.tamu.edu/publications-reports/tr/0110.pdf}.

\bibitem[Veeser and Zanotti(2018{\natexlab{a}})]{VeeZa:18}
A.~Veeser and P.~Zanotti.
\newblock Quasi-optimal nonconforming methods for symmetric elliptic problems.
  {I}---{A}bstract theory.
\newblock \emph{SIAM J. Numer. Anal.}, 56\penalty0 (3):\penalty0 1621--1642,
  2018{\natexlab{a}}.

\bibitem[Veeser and Zanotti(2018{\natexlab{b}})]{VeeZa:18b}
A.~Veeser and P.~Zanotti.
\newblock Quasi-optimal nonconforming methods for symmetric elliptic problems.
  {III}---{D}iscontinuous {G}alerkin and other interior penalty methods.
\newblock \emph{SIAM J. Numer. Anal.}, 56\penalty0 (5):\penalty0 2871--2894,
  2018{\natexlab{b}}.

\end{thebibliography}

\end{document}